\documentclass[times,sort&compress,3p]{elsarticle}
\journal{Journal of Multivariate Analysis}
\usepackage[labelfont=bf]{caption}

\usepackage{pdflscape}
\usepackage{longtable}
\usepackage{graphicx}
\usepackage{multicol}
\usepackage{csquotes}
\usepackage{dsfont}
\usepackage{pdflscape}
\usepackage{morefloats}
\usepackage{nccfoots}
\usepackage{enumitem} 
\usepackage{booktabs}
\usepackage{inputenc}
\usepackage{hyperref}
\usepackage{lscape}
\usepackage{multirow}
\usepackage{geometry,blindtext}

\usepackage{amsmath,amsfonts,amssymb,amsthm,booktabs,color,epsfig,graphicx,hyperref,url}

\theoremstyle{plain}% Theorem-like structures provided by amsthm.sty
\newtheorem{theorem}{Theorem}

\newtheorem{proposition}{Proposition}
\newtheorem{lemma}{Lemma}

\theoremstyle{definition}
\newtheorem{definition}{Definition}
\newtheorem{remark}{Remark}
\newtheorem{example}{Example}

\renewenvironment{proof}[1][\proofname]{{\noindent\bfseries #1. }}{\qed}

\newcommand{\E}{{\operatorname E}}
\newcommand{\Var}{{\operatorname{Var}}}
\newcommand{\Cov}{{\operatorname{Cov}}}
\newcommand{\Cor}{{\operatorname{Corr}}}
\newcommand{\tr}{{\operatorname{tr}}}
\newcommand{\OPD}{{\operatorname{OPD}}}
\newcommand{\arp}{\operatorname{AR}}

\newcommand{\bX}{\mathbf{X}}
\newcommand{\bY}{\mathbf{Y}}
\newcommand{\R}{\mathbb{R}}
\newcommand{\Z}{\mathbb{Z}}
\newcommand{\A}{\mathcal{A}}
\newcommand{\B}{\mathcal{B}}
\newcommand{\F}{\mathcal{F}}

\newcommand{\pr}[1]{\left(#1\right)}
\newcommand{\Prob}{\operatorname{Pr}}

\makeatletter
\newcommand*{\defeq}{\mathrel{\rlap{%
                     \raisebox{0.3ex}{$\m@th\cdot$}}%
                     \raisebox{-0.3ex}{$\m@th\cdot$}}%
                     =}
\makeatother

\makeatletter
\newcommand*{\eqdef}{=\mathrel{\rlap{%
                     \raisebox{0.3ex}{$\m@th\cdot$}}%
                     \raisebox{-0.3ex}{$\m@th\cdot$}}%
                     }
\makeatother

\usepackage{listings}
\begin{document}

\begin{frontmatter}

 \title{Ordinal pattern dependence as a multivariate dependence measure}
 
 \author[1]{Annika Betken}
 \author[2]{Herold Dehling}
 \author[3]{Ines N\"ußgen} 
 \author[4]{Alexander Schnurr\corref{mycorrespondingauthor}}
% \ead[url]{home page}

\cortext[mycorrespondingauthor]{Corresponding author. Email address: \url{schnurr@mathematik.uni-siegen.de}}

\address[1]{Faculty of Electrical Engineering, Mathematics and Computer Science, University of Twente,  7500 AE Enschede, The Netherlands}
\address[2]{Faculty of Mathematics, Ruhr-University Bochum, 44780 Bochum, Germany}
\address[3]{Department of Mathematics,  Siegen University, 57072 Siegen, Germany}
\address[4]{Department of Mathematics,  Siegen University, 57072 Siegen, Germany}

\begin{abstract}
In this article, we show that the recently introduced ordinal pattern dependence fits into the axiomatic framework of general multivariate dependence measures, i.e., 
 measures of dependence between two multivariate random objects.
Furthermore, we consider   multivariate generalizations of established univariate dependence measures like  Kendall's $\tau$, Spearman's $\rho$ and Pearson's correlation coefficient.
Among these, only  multivariate Kendall's $\tau$
proves to take the dynamical dependence of random vectors 
stemming from multidimensional time series
into account. Consequently, the article focuses on a comparison of ordinal pattern dependence and multivariate  Kendall's $\tau$  in this context.
To this end,  limit theorems for multivariate Kendall's $\tau$
are established under the assumption of near-epoch dependent data-generating time series. We analyze how
ordinal pattern dependence compares to multivariate Kendall's $\tau$ and Pearson's correlation coefficient on theoretical grounds.
 Additionally, a simulation study illustrates differences in the kind of dependencies that are revealed by multivariate Kendall's $\tau$ and ordinal pattern dependence.
\end{abstract}


\end{frontmatter}

\section{Introduction\label{sec:1}}

Recently, various attempts have been made to generalize classical dependence measures for one-dimensional random variables (like Pearson's correlation coefficient, Kendall's $\tau$, Spearman's $\rho$) to a multivariate framework. The aim of these  is to describe the degree of dependence between two random vectors with a single number.  This has to be separated from the branch of research where the dependence \emph{within} one vector is described by a single number (see  \cite{quessy:2013, Schmid_et_al:2010} and the references therein). 

Roughly speaking, one can separate the following two approaches: (I) In a first step, the main properties which classical dependence measures between two random variables display,   are extracted. In a second step, multivariate analogues of the dependence measures which satisfy canonical generalizations of these properties in a multivariate framework, are defined. However, often  a canonical interpretation of these measures is  not at hand. (II) Given two time series, one wants to describe their co-movement. 

Along these lines, the definition of ordinal pattern dependence (see \cite{schnurr:2014}) follows the latter approach. Originally, axiomatic systems are disregarded by the notion of ordinal pattern dependence, which  is naturally interpreted as the degree of co-monotonic behavior of two time series.
Against the background of this approach, limit theorems have been proved in the time series setting (see \cite{schnurr:dehling:2017} for the SRD case and \cite{nuessgen:schnurr:2021} for the LRD case). 

Both approaches in defining multivariate dependence measures have proved to be useful, but by now, they have been analyzed separately. In the present paper, we close the gap between the two. 
To this end,
we recall the definition of ordinal pattern dependence in the subsequent section and show that it is a multivariate dependence measure according to the definition introduced in \cite{grothe:2014}. In Section \ref{sec:iid}, we establish consistency and asymptotic normality for estimators of ordinal pattern dependence in the framework of i.i.d. random vectors. Section \ref{sec:othermeasures} deals with multivariate extensions of well-established univariate dependence measures. It turns out that multivariate Kendall's $\tau$ is the only one among these that captures the dynamical dependence between random vectors. Starting with approach (I), we prove limit theorems for an estimator of multivariate Kendall's $\tau$ in the time series context. In the last section, the different measures are compared from a theoretical point-of-view as well as by simulation studies. 

\section{Ordinal pattern dependence as a measure of multivariate dependence}
\label{sec:axiomatic}

%%%%%%%%%%%%%%%
%%%%%%%%%%%%%%%%%

If $(X_i,Y_i)$, $i\geq 1$, denotes a stationary, bivariate process, we define, for any integers $i,h\geq 1$,  the random vectors of consecutive observations
\begin{eqnarray*}
  \bX_i^{(h)} &\defeq & (X_i,\ldots,X_{i+h}), \\
  \bY_i^{(h)} &\defeq  & (Y_i,\ldots,Y_{i+h}).
\end{eqnarray*}

The goal of this paper is to consider the concept of ordinal pattern dependence as a multivariate measure of dependence between the random vectors $\bX_i^{(h)}$ and $\bY_i^{(h)}$ stemming from a stationary, bivariate process $(X_i,Y_i)$, $i\geq 1$, with continuous marginal distributions, and to compare it to established measures of dependence. 
Note that, by stationarity of the underlying process, the joint distribution of the vector $(\bX_i^{(h)},\bY_i^{(h)})$ does not depend on $i$. We will thus use the  symbol $(\bX^{(h)},\bY^{(h)})$ for a generic random vector with the same joint distribution as any of the $(\bX_i^{(h)},\bY_i^{(h)})$ and we write
 $\bX^{(h)} = (X_1,\ldots,X_{1+h})$,
 $\bY^{(h)} = (Y_1,\ldots,Y_{1+h})$.
Furthermore, note that it is common to count the number of increments $h$ rather than the length of the vector, since ordinal patterns can be calculated by exclusively considering the increments of the time series. Moreover, here, and in the following, we consider vectors as column vectors. However, for the sake of readability and notational convenience, we omit the notation $\top$ indicating the transpose of vectors.

\subsection{Ordinal pattern dependence}

For $h\in\mathbb{N}$ let $S_{h}$ denote the set of permutations of $\{0, \ldots, h\}$, which we write as $(h+1)$-tuples containing each of the numbers $0, \ldots, h$ exactly once.
The ordinal pattern of order $h$  refers to the permutation
\begin{align*}
\Pi(x_0, \ldots, x_h)=(\pi_0,\ldots, \pi_h)\in S_h
\end{align*}
 which satisfies
$
x_{\pi_0}\geq \cdots\geq x_{\pi_h}
$
(see \cite{bandt:pompe:2002, bandt:Shiha:2007}). In the present paper, we only consider continuous marginals.
Allowing for non-continuous marginals would require the additional restriction 
$\pi_{j-1}>\pi_j$ if $x_{\pi_{j-1}}=x_{\pi_j}$ for $j\in \{1,...,h\}$ (see \cite{sinn:keller:2011}).

\begin{definition}
We define the ordinal pattern dependence between two random vectors  $\bX^{(h)}= (X_1, \ldots, X_{h+1})$ and $\bY^{(h)}= (Y_1, \ldots, Y_{h+1})$ by
\begin{align} \label{posdep}
\OPD_h(\bX^{(h)}, \bY^{(h)}) =\frac{\Prob\left(\Pi\pr{\bX^{(h)}}=\Pi\pr{\bY^{(h)}}\right)-\sum_{\pi\in S_{h}}\Prob\left(\Pi(\bX^{(h)})=\pi\right)\Prob\left(\Pi(\bY^{(h)})=\pi\right)}{1-\sum_{\pi\in S_{h}}\Prob\left(\Pi(\bX^{(h)})=\pi\right)\Prob\left(\Pi(\bY^{(h)})=\pi\right)}.
\end{align}
\end{definition}
This definition of ordinal pattern dependence only takes positive dependence into account.
 Negative dependence can be included by analyzing the co-movement of $X$ and $-Y=(-Y_i)_{i\in\mathbb{N}}$. 
Typically, one is interested in measuring either positive \emph{or} negative dependence. If one wants to consider both dependencies at the same time, a consideration of the quantity
\[
\OPD_h(\bX, \bY)^+-\OPD_h(\bX, -\bY)^+,
\]
where $a^+:=\max\{ a,0 \}$ for every $a\in\mathbb{R}$, seems natural. In order to keep things less technical, we only consider the simpler measure \eqref{posdep}. For recent developments in the theory of ordinal patterns, see \cite{bandt:2020} and  \cite{mohr:2020},  and  for a related approach to analyze dependence between dynamical systems, see \cite{echegoyen_et_al:2019}. 

 \subsection{Axiomatic definition of multivariate dependence measures}

With the following definition,
\cite{grothe:2014} establish an axiomatic theory for multivariate dependence measures between $n$-dimensional random vectors.  This has been strongly inspired by the axiomatic framework of \cite{Schmid_et_al:2010}, who follow a copula-based approach to define and analyze multivariate dependence measures within one vector.

\begin{definition}\label{def:dep_measure}
Let $L_0$ denote the space of random vectors with values in $\mathbb{R}^n$ on the common probability space $(\Omega, \F, \Prob)$. We call a function $\mu: L_0 \times L_0 \longrightarrow\mathbb{R}$  an $n$-dimensional measure of dependence if
\begin{enumerate}
\item it takes values in $[-1,1]$; 
\item  it is invariant with respect to simultaneous permutations of the components within two random vectors $\mathbf{X}$ and $\mathbf{Y}$;
\item  it is invariant with respect to monotonically increasing transformations of the components of the  two random vectors $\mathbf{X}$ and $\mathbf{Y}$; 
\item  it is zero for two independent random vectors $\mathbf{X}$ and $\mathbf{Y}$;
\item  it respects concordance ordering, i.e., for two pairs of random vectors  
 $\mathbf{X}$, $\mathbf{Y}$ and $\mathbf{X}^*$, $\mathbf{Y}^*$, it holds that
\[
  \binom{\mathbf{X}} {\mathbf{Y}} \preccurlyeq_C  \binom{\mathbf{X}^*} {\mathbf{Y}^*} \Rightarrow \mu(\mathbf{X},\mathbf{Y})\leq \mu(\mathbf{X}^*,\mathbf{Y}^*).
\]
Here, $\preccurlyeq_C $ denotes concordance ordering, i.e., 
\[
  \binom{\mathbf{X}} {\mathbf{Y}} \preccurlyeq_C  \binom{\mathbf{X}^*} {\mathbf{Y}^*} \ \text{if and only if} \
F_{\binom{\mathbf{X}}{\mathbf{Y}}} \leq F_{\binom{\mathbf{X}^*}{\mathbf{Y}^*}} \text{ and }
\bar{F}_{\binom{\mathbf{X}}{\mathbf{Y}}} \leq \bar{F}_{\binom{\mathbf{X}^*}{\mathbf{Y}^*}},
\]
where $\leq$ is meant pointwise and $\bar{F}$ denotes the survival function.
\end{enumerate}
\end{definition}

\begin{theorem}\label{thm:opd_dep_measure}
The ordinal pattern dependence $\OPD_h$ is an $h+1$-dimensional measure of dependence.
\end{theorem}

The proof, which is a bit involved and makes use of mulitvariate distribution functions and survival functions, has been postponed to Section \ref{sec:proofs}. 

\section{Limit Theorems for Ordinal Pattern Dependence of i.i.d. Vectors}
\label{sec:iid}

In Section \ref{sec:comparison}, we compare ordinal pattern dependence to other concepts of multivariate dependence. These have been introduced and used for sequences of independent random vectors. In contrast to this, the definition of ordinal pattern dependence applies to random vectors stemming from multivariate time series.
Nonetheless, ordinal pattern dependence can as well 
be applied to independent random vectors. 
Limit theorems that provide the asymptotic distribution of ordinal pattern dependence in this setting have not yet been established.
We close this gap by the following considerations:

Let $(\bX_i,\bY_i)$, $i\geq 1$,  be independent copies of $(\bX, \bY)$, and define
\begin{eqnarray*}
 \hat{q}_{\bX, \pi, n}:=\frac{1}{n}\sum\limits_{i=1}^{n}1_{\left\{\Pi(\bX_{i})=\pi\right\}}, \
  \hat{q}_{\bY,  \pi, n}:= \frac{1}{n}\sum\limits_{i=1}^{n}1_{\left\{\Pi(\bY_{i})=\pi\right\}}, \
 \hat{q}_{(\bX,\bY),  n}:= \frac{1}{n}\sum\limits_{i=1}^{n}1_{\left\{\Pi(\bX_{i})=\Pi(\bY_{i})\right\}}, 
\end{eqnarray*}
as well as the corresponding probabilities
\begin{eqnarray*}
q_{\bX^{}, \pi}\defeq \Prob\left(\Pi(\bX^{})=\pi\right),\
q_{\bY^{}, \pi} \defeq \Prob\left(\Pi(\bY^{ })=\pi\right),\
q_{(\bX^{ },\bY^{ })} \defeq  \Prob \left(\Pi(\bX^{ })=\Pi(\bY^{ })\right).
\end{eqnarray*}
 
According to the law of large numbers  $\hat{q}_{\bX,\pi,n}$,
$\hat{q}_{\bY,\pi,n}$, and  $\hat{q}_{(\bX,\bY),n}$  are strongly consistent estimators for these probabilities.

\begin{proposition}
Let  $(\bX_i,\bY_i)$, $i\geq 1$,   be independent  copies of $(\bX^{ }, \bY^{ })$.  
Then, as $n\rightarrow\infty$,
\begin{eqnarray*}
 \hat{q}_{\bX, \pi, n}\longrightarrow q_{\bX^{ },\pi}, \
  \hat{q}_{\bY , \pi ,n}\longrightarrow q_{\bY^{ },\pi}, \
 \hat{q}_{(\bX,\bY),n}\longrightarrow q_{(\bX^{ },\bY^{ })}
\end{eqnarray*}
almost surely. 
\end{proposition}

The following theorem establishes asymptotic normality of ordinal pattern dependence of i.i.d. random vectors.
For this, we introduce the following notation:
\begin{eqnarray*}
\hat{q}_{\bX,n}:= \left( \hat{q}_{\bX,\pi,n}  \right)_{\pi\in \mathcal{S}_{h+1}}, \
\hat{q}_{\bY,n}:= \left( \hat{q}_{\bY,\pi,n}  \right)_{\pi\in \mathcal{S}_{h+1}}, \
q_{\bX}:= \left(  q_{\bX,\pi} \right)_{\pi\in \mathcal{S}_{h+1}}, \
q_{\bY}:= \left(  q_{\bY,\pi} \right)_{\pi\in \mathcal{S}_{h+1}}.
\end{eqnarray*} 

\begin{theorem}\label{th:opd-iid-an}
Let  $(\bX_i,\bY_i)$, $i\geq 1$,  be independent copies of $(\bX^{ }, \bY^{ })$.  
Then, as $n\rightarrow\infty$,
\[
\sqrt{n}\left( \frac{ \hat{q}_{(\bX,\bY) ,n}-\sum_{\pi\in S_{h}} \hat{q}_{\bX, \pi, n} \hat{q}_{\bY, \pi, n}}{1-\sum_{\pi\in S_{h}} \hat{q}_{\bX^{ }, \pi, n} \hat{q}_{\bY^{ }, \pi, n}}
- 
 \frac{q_{(\bX^{ },\bY^{ })}-\sum_{\pi\in S_{h}} q_{\bX^{ }, \pi} q_{\bY^{ }, \pi}}{1-\sum_{\pi\in S_{h}} q_{\bX^{ }, \pi} q_{\bY^{ }, \pi}}
\right)
\stackrel{\mathcal{D}}{\longrightarrow} N(0,\sigma^2),
\]
where the limit variance $\sigma^2$ is given by
\[
  \sigma^2= \nabla f (q_{(\bX,\bY)}, q_{\bX},q_{\bY})\, \Sigma\, \big(\nabla f (q_{(\bX,\bY)}, q_{\bX},q_{\bY})\big)^\top.
\]
Here, the matrix $\Sigma$ is defined as in Proposition~\ref{prop:opd-iid-an} (see below), and $\nabla f$ is the gradient of the function
\newline $f:\R\times \R^{(h+1)!}\times \R^{(h+1)!}\rightarrow \R$, defined by
\[
  f(u,v,w)=\frac{u-v^\top\cdot w}{1-v^\top\cdot w}.
\]
\end{theorem}

The proof of Theorem~\ref{th:opd-iid-an} is based on the following proposition, which establishes the joint asymptotic normality  of  $\hat{q}_{\bX,\pi,n}$, $\hat{q}_{\bY,\pi,n}$, and  $\hat{q}_{(\bX,\bY),n}$. 

\begin{proposition}\label{prop:opd-iid-an}
Under the same assumptions as in Theorem~\ref{th:opd-iid-an}, we have
\begin{equation*}
\sqrt{n} \left( 
\begin{array}{c} 
 \hat{q}_{(\bX,\bY),n}-q_{(\bX,\bY)} \\
  \hat{q}_{\bX,n} -q_{\bX} \\
   \hat{q}_{\bY,n} -q_{\bY} 
\end{array} 
\right) \stackrel{\mathcal{D}}{\longrightarrow} N(0,\Sigma),
\end{equation*} 
where $\Sigma$  is the symmetric $(2 (h+1)!+1)\times (2 (h+1)! +1)$  matrix  
\[
 \Sigma=\left(
  \begin{array}{ccc}
   \Sigma_{11} & \Sigma_{12} & \Sigma_{13} \\
    \Sigma_{21} & \Sigma_{22} & \Sigma_{23} \\
     \Sigma_{31} & \Sigma_{32} & \Sigma_{33}
  \end{array} 
   \right)
\]
with 
\begin{align*}
&\Sigma_{11}=q_{(\bX,\bY)}(1-q_{(\bX,\bY)}) \in \R, \\
&\Sigma_{12} =(\sigma_\pi(1,2) )_{\pi\in \mathcal{S}_{h}}\in \R^{1\times (h+1)!},  \ 
 &&\sigma_{\pi}(1,2)=\Prob(\Pi(X)=\Pi(Y)=\pi)-q_{\bX,\bY}\cdot q_{\bX,\pi},
\\
&\Sigma_{13}=(\sigma_\pi(1,3))_{\pi\in \mathcal{S}_{h}}\in \R^{1\times (h+1)!},  \
&& \sigma_{\pi}(1,3)=\Prob(\Pi(X)=\Pi(Y)=\pi)-q_{\bX,\bY}\cdot q_{\bY,\pi},\\
&\Sigma_{22}=(\sigma_{\pi,\pi^\prime}(2,2) )_{\pi,\pi^\prime\in \mathcal{S}_{h}}\in \R^{(h+1)!\times (h+1)!},  \
&&\sigma_{\pi,\pi^\prime}(2,2)=
\left\{ 
\begin{array}{ll} 
 -q_{\bX,\pi}q_{\bX,\pi^\prime} & \mbox{ if } \pi \neq \pi^\prime \\
 q_{\bX,\pi}(1-q_{\bX,\pi}) & \mbox{ if } \pi = \pi^\prime,
\end{array} 
\right.\\
 &\Sigma_{23}=(\sigma_{\pi,\pi^\prime}(2,3) )_{\pi,\pi^\prime\in \mathcal{S}_{h}}\in \R^{(h+1)!\times (h+1)!}, \
&&\sigma_{\pi,\pi^\prime}(2,3)= \Prob(\Pi(\bX)=\pi, \Pi(\bY)=\pi^\prime) - q_{\bX,\pi} q_{\bY,\pi^\prime}, \\
&\Sigma_{33}=(\sigma_{\pi,\pi^\prime}(3, 3) )_{\pi,\pi^\prime\in \mathcal{S}_{h}}\in \R^{(h+1)!\times (h+1)!},  \
&&\sigma_{\pi,\pi^\prime}(3, 3)=
\left\{ 
\begin{array}{ll} 
 -q_{\bY,\pi}q_{\bY,\pi^\prime} & \mbox{ if } \pi \neq \pi^\prime \\
 q_{\bY,\pi}(1-q_{\bY,\pi}) & \mbox{ if } \pi = \pi^\prime.
\end{array} 
\right.
\end{align*} 
Due to symmetry of $\Sigma$, the remaining blocks are defined by $\Sigma_{21}=\Sigma_{12}^\top$, $\Sigma_{31}=\Sigma_{13}^\top$,  $\Sigma_{32}=\Sigma_{23}^\top$.  
\end{proposition}
\begin{proof}
The proof follows directly  from the multivariate central limit theorem applied to the partial sums of the \newline $(2(h+1)!+1)$-dimensional i.i.d. random vectors
\[
 \xi_i= \left( 1_{\{\Pi(\bX_i)=\Pi(\bY_i)\}},(1_{\{\Pi(\bX_i)=\pi\}})_{\pi\in \mathcal{S}_{h}}, 
  (1_{\{ \Pi(\bY_i)=\pi \}})_{\pi\in \mathcal{S}_{h}}\right).
\]
The limit covariance matrix is the covariance matrix of $\xi_1$, which is given by the formulae stated in the formulation of this proposition.
\end{proof}

\begin{proof}[Proof of Theorem~\ref{th:opd-iid-an}] We apply the delta method to the function $f$, defined in the formulation of the theorem, together with the multivariate CLT established in Proposition~\ref{prop:opd-iid-an}. In this way, we obtain
\begin{eqnarray*}
   \sqrt{n}\left( f(\hat{q}_{(\bX,\bY),n},\hat{q}_{\bX,n},\hat{q}_{\bY,n}) - f(q_{(\bX^{(h)},\bY^{(h)})},q_{\bX^{(h)}},q_{\bY^{(h)}})  \right) 
  \stackrel{\mathcal{D}}{\longrightarrow} 
  N(0,\nabla f(q_{(\bX^{(h)},\bY^{(h)})},q_{\bX^{(h)}},q_{\bY^{(h)}}) \cdot \Sigma \cdot( \nabla f(q_{(\bX^{(h)},\bY^{(h)})},q_{\bX^{(h)}},q_{\bY^{(h)}}))^\top).
\end{eqnarray*}
This proves the statement of the theorem.
\end{proof}

 \section{Ordinal pattern dependence in contrast to multivariate Kendall's $\tau$}
\label{sec:othermeasures}

In this article, we are explicitly studying the dependence between random vectors stemming from stationary time series. 
In this regard, the main drawback of  univariate dependence measures  is that these do not incorporate cross-dependencies which characterize the \textit{dynamical} dependence between two random vectors. Univariate dependence measures focus on the dependence between $X_i$ and $Y_i$, i.e., on the dependence   at the same point in time.
In contrast, ordinal pattern dependence captures the dynamics of time series.

In the following, we study two multivariate generalizations of univariate dependence measures, namely the multivariate extension of Pearson's correlation coefficient, established in \cite{puccetti:2019}, and multivariate Kendall's $\tau$ as introduced in \cite{grothe:2014}.

\begin{definition}
For two $(h+1)$-dimensional random vectors $\bX^{ }, \bY^{ }\in L_2$ with invertible covariance matrices $\Sigma_{\bX^{ }}$ and $\Sigma_{\bY^{ }}$ and cross-covariance matrix $\Sigma_{\bX^{ },\bY^{ }}$, we define Pearson's correlation coefficient by
\begin{align*}
\rho\left(\bX^{ }, \bY^{ }\right)\defeq\frac{\tr\left(\Sigma_{\bX^{ },\bY^{ }}\right)}{\tr\left(\left(\Sigma_{\bX^{ }}\Sigma_{\bY^{ }}\right)^{1/2}\right)},
\end{align*}
where $A^{1/2}$ is the principal square root of the matrix $A$, such that $A^{1/2}A^{1/2}=A$.
\end{definition}
 
For the multivariate generalization of Pearson's  correlation coefficient,  we obtain
\begin{align*}
\rho\left( \bX^{ }, \bY^{ }\right)=\frac{\tr\left(\Sigma_{ \bX^{ }, \bY^{ }}\right)}{\tr\left(\left(\Sigma_{ \bX^{ }}\Sigma_{ \bY^{ }}\right)^{1/2}\right)}=\frac{\Cov(X_1,Y_1)+\ldots+\Cov(X_{1+h},Y_{1+h})}{\tr\left(\left(\Sigma_{ \bX^{ }}\Sigma_{ \bY^{ }}\right)^{1/2}\right)}.
\end{align*}
As a result, the cross-correlations have no impact on the value of  Pearson's  correlation coefficient.
%If we additionally assume that  $\Sigma_{ \bX^{(h)}}=\Sigma_{ \bY^{(h)}}$, it follows that
%\begin{align*}
%\rho\left( \bX^{(h)}, \bY^{(h)}\right)=\Cov(X_1,Y_1).
%\end{align*}
Therefore, the multivariate Pearson's  correlation coefficient does not seem to be  appropriate for our approach. 
The same holds true for generalizations of Spearman's $\rho$ due to the close relationship between these concepts. 
We hence focus on the multivariate generalization of Kendall's $\tau$:

\subsection{Multivariate Kendall's $\tau$}
The definition of multivariate Kendall's $\tau$ that we consider in this section 
is taken from \cite{grothe:2014}.  In that paper, the authors investigated the dependence between two multivariate random vectors. Therefore, for our purposes, it is appropriate to use it in the time series context. For a multivariate generalization of Kendall's $\tau$ \textit{within} one random vector see \cite{quessy:2013}.
\begin{definition}
For two $h+1$-dimensional random vectors $\bX^{ }, \bY^{ }$, we define Kendall's $\tau$ by
\begin{align*}
\tau(\bX^{ },\bY^{ })\defeq
\Cor\left(\mathbf{1}_{\left\{ \bX^{ }\leq \tilde{\bX}^{ }\right\}},\mathbf{1}_{\left\{ \bY^{ }\leq \tilde{\bY}^{ }\right\}} \right)
,
\end{align*}
where $\left(\tilde{\bX}^{ },\tilde{\bY}^{ }\right)$ is an independent copy of $(\bX^{ },\bY^{ })$.
\end{definition}

The following lemma establishes a representation of  multivariate Kendall's $\tau$ for Gaussian processes
in terms  of the probabilities $p_{\bX^{ }}$ and $p_{\bY^{ }}$ that enter in our definition of ordinal pattern dependence.
\begin{lemma} \label{lemma1}
Let $(X_i,Y_i)$, $i\geq 1$,  denote a stationary mean zero Gaussian process and let   $\bX^{(h)} = (X_1,\ldots,X_{1+h})$ and
  $\bY^{(h)} = (Y_1,\ldots,Y_{1+h})$. Then, we have
\begin{align*}
\tau\left( \bX^{(h)}, \bY^{(h)}\right)=
\frac{\Prob\left(X_1\leq 0,\ldots,X_{1+h}\leq 0,Y_1\leq 0,\ldots,Y_{1+h}\leq 0\right)-\tilde{p}_{\bX^{(h)}}\tilde{p}_{\bY^{(h)}}}{\sqrt{\tilde{p}_{\bX^{(h)}}\left(1-\tilde{p}_{\bX^{(h)}}\right)\tilde{p}_{\bY^{(h)}}\left(1-\tilde{p}_{\bY^{(h)}}\right)}},
\end{align*}
where $\tilde{p}_{\bX^{(h)}}=\Prob\left(X_1\leq 0,\ldots,X_{1+h}\leq 0\right)$ and $\tilde{p}_{\bY^{(h)}}=\Prob\left(Y_1\leq 0,\ldots,Y_{1+h}\leq 0\right)$. 
\end{lemma}

\begin{proof}
Let $\left(\tilde{\bX}^{(h)},\tilde{\bY}^{(h)}\right)$ be an independent copy of $(\bX^{(h)},\bY^{(h)})$ with $\tilde{\bX}^{(h)}=(\tilde{X}_1,\ldots,\tilde{X}_{1+h})$ and $\tilde{\bY}^{(h)}=(\tilde{Y}_1, \ldots,\tilde{Y}_{1+h})$.
It then  holds that
\begin{align*}
\tau\left( \bX^{(h)}, \bY^{(h)}\right)
&=\Cor\left(\mathbf{1}_{\left\{ \bX^{(h)}\leq \tilde{\bX}^{(h)}\right\}},\mathbf{1}_{\left\{ \bY^{(h)}\leq \tilde{\bY}^{(h)}\right\}} \right)\\
&=\Cor\left(\mathbf{1}_{\left\{ \bX^{(h)}- \tilde{\bX}^{(h)}\leq 0\right\}},\mathbf{1}_{\left\{ \bY^{(h)}-\tilde{\bY}^{(h)}\leq 0\right\}} \right)\\
&=\frac{\Prob\left(X_1-\tilde{X}_1\leq 0,\ldots,X_{1+h}-\tilde{X}_{1+h}\leq 0,Y_1-\tilde{Y}_1\leq 0,\ldots,Y_{1+h}-\tilde{Y}_{1+h}\leq 0\right)-p_{\bX^{(h)}}p_{\bY^{(h)}}}{\sqrt{p_{\bX^{(h)}}\left(1-p_{\bX^{(h)}}\right)p_{\bY^{(h)}}\left(1-p_{\bY^{(h)}}\right)}}
\end{align*}
with
$p_{\bX^{(h)}}=\Prob\left(X_1-\tilde{X}_1\leq 0,\ldots,X_{1+h}-\tilde{X}_{1+h}\leq 0\right)$
and
$p_{\bY^{(h)}}=\Prob\left(Y_1-\tilde{Y}_1\leq 0,\ldots,Y_{1+h}-\tilde{Y}_{1+h}\leq 0\right)$.
Note that for independent centered Gaussian processes,
\begin{align*}
\left( \bX^{(h)}-\tilde{\bX}^{(h)}, \bY^{(h)}-\tilde{\bY}^{(h)}\right)\overset{D}{=}\sqrt{2}\left( \bX^{(h)}, \bY^{(h)}\right).
\end{align*}
This explicitly implies that the cross-correlations within  $\left( \bX^{(h)}-\tilde{\bX}^{(h)}, \bY^{(h)}-\tilde{\bY}^{(h)}\right)$ equal those within $\left( \bX^{(h)}, \bY^{(h)}\right)$. Therefore, we have $p_{\bX^{(h)}}=\tilde{p}_{\bX^{(h)}}$ and 
\begin{align*}
&\hspace*{-10mm}\frac{\Prob\left(X_1-\tilde{X}_1\leq 0,\ldots,X_{1+h}-\tilde{X}_{1+h}\leq 0,Y_1-\tilde{Y}_1\leq 0,\ldots,Y_{1+h}-\tilde{Y}_{1+h}\leq 0\right)-p_{\bX^{(h)}}p_{\bY^{(h)}}}{\sqrt{p_{\bX^{(h)}}\left(1-p_{\bX^{(h)}}\right)p_{\bY^{(h)}}\left(1-p_{\bY^{(h)}}\right)}}\\
=&\frac{\Prob\left(X_1\leq 0,\ldots,X_{1+h}\leq 0,Y_1\leq 0,\ldots,Y_{1+h}\leq 0\right)-\tilde{p}_{\bX^{(h)}}\tilde{p}_{\bY^{(h)}}}{\sqrt{\tilde{p}_{\bX^{(h)}}\left(1-\tilde{p}_{\bX^{(h)}}\right)\tilde{p}_{\bY^{(h)}}\left(1-\tilde{p}_{\bY^{(h)}}\right)}}.
\end{align*} 
\end{proof}

Although for $h\geq 2$ we cannot derive an analytic expression for 
\begin{align*}
\Prob\left(X_1\leq 0,\ldots,X_{1+h}\leq 0\right), \ \Prob\left(Y_1\leq 0,\ldots,Y_{1+h}\leq 0\right)\\
\intertext{or}
\Prob\left(X_1\leq 0,\ldots,X_{1+h}\leq 0,Y_1\leq 0,\ldots,Y_{1+h}\leq 0\right),
\end{align*}
we know that these orthant probabilities of a multivariate Gaussian distribution are determined by the entries of the correlation matrices and by the entries of the cross-correlation matrix of $\bX^{(h)}$ and $\bY^{(h)}$. 
In contrast to multivariate Pearson's  correlation coefficient, multivariate Kendall's $\tau$
constitutes a multivariate dependence measure that takes the dynamical dependence of data stemming from time series into account.

\subsection{Estimation of multivariate Kendall's $\tau$}

\cite{grothe:2014} consider an estimator for  multivariate Kendall's $\tau$ based on independent vectors $(\bX_i,\bY_i)$, $1\leq i\leq n$. In our setup, we will define an empirical version of Kendall's $\tau$ based on the dependent vectors 
$(\bX_i^{(h)},\bY_i^{(h)})$, $1\leq i\leq n$. For this, we will follow the ideas of \cite{dehling:2017}, who considered estimation of the classical univariate Kendall's $\tau$ for bivariate time series under some mild dependence condition.

Given an independent copy $(\tilde{\bX}^{(h)},\tilde{\bY}^{(h)})$  of the vector $(\bX^{(h)},\bY^{(h)})$, we have
\begin{eqnarray*}
 \tau(\bX^{(h)},\bY^{(h)})
 &=& \frac{ p_{(\bX^{(h)},\bY^{(h)})} -p_{\bX^{(h)}} p_{\bY^{(h)}}}{ 
 \sqrt{  p_{\bX^{(h)}} (1- p_{\bX^{(h)}}  ) p_{\bY^{(h)}} (1- p_{\bY^{(h)}})}} \\
 &=& \psi( p_{\bX^{(h)}}, p_{\bY^{(h)}},p_{(\bX^{(h)},\bY^{(h)})} ),
 \end{eqnarray*}
where $p_{(\bX^{(h)},\bY^{(h)})}:=\Prob(\bX^{(h)} \leq \tilde{\bX}^{(h)},  \bY^{(h)}\leq \tilde{\bY}^{h})$, 
$p_{\bX^{(h)}} := \Prob(\bX^{(h)} \leq \tilde{\bX}^{(h)})$, $p_{\bY^{(h)}} :=\Prob(\bY^{(h)} \leq \tilde{\bY}^{(h)})$, and where $\psi:\R^3\rightarrow \R$ is defined by
\begin{align}\label{eq:psi}
  \psi(x,y,z):= \frac{z-x\, y}{\sqrt{x(1-x)y(1-y)} }. 
\end{align}

The probabilities $p_{\bX^{(h)}}$, $p_{\bY^{(h)}}$, and $p_{(\bX^{(h)},\bY^{(h)})}$ can be estimated by their sample analogues defined by
\begin{align*}
 \hat{p}_{\bX^{(h)},n}:= \frac{1}{n(n-1)} 
 \sum_{1\leq i \neq j\leq n} 1_{\{ \bX_i^{(h)} \leq \bX_j^{(h)} \}}, \
  &\hat{p}_{\bY^{(h)},n}:= \frac{1}{n(n-1)} 
  \sum_{1\leq i \neq j\leq n} 1_{\{ \bY_i^{(h)} \leq \bY_j^{(h)}  \}}, \\
 \hat{p}_{(\bX^{(h)},\bY^{(h)}),n}:= \frac{1}{n(n-1)}  &
 \sum_{1\leq i \neq j\leq n} 1_{\{ \bX_i^{(h)} \leq \bX_j^{(h)}, \bY_i^{(h)} \leq \bY_j^{(h)}  \}}, 
\end{align*}
where   $\bX_i^{(h)} = (X_i,\ldots,X_{i+h})$ and
  $\bY_i^{(h)} = (Y_i,\ldots,Y_{i+h})$.
The plug-in estimator for Kendall's $\tau$ is then given by 
\[
  \hat{\tau}_n(\bX^{(h)},\bY^{(h)}):=\psi(\hat{p}_{\bX^{(h)},n},  \hat{p}_{\bY^{(h)},n},  \hat{p}_{(\bX^{(h)},\bY^{(h)}),n}).
\]
In what follows, we will derive the joint limit distribution of the random vector $(\hat{p}_{\bX^{(h)},n},  \hat{p}_{\bY^{(h)},n},  \hat{p}_{(\bX^{(h)},\bY^{(h)}),n})$ and  the limit distribution of $\hat{\tau}_n(\bX^{(h)},\bY^{(h)})$ by the delta method. 
For this, observe that  $\hat{p}_{\bX^{(h)},n}$,  $\hat{p}_{\bY^{(h)},n}$, and   $\hat{p}_{(\bX^{(h)},\bY^{(h)}),n}$ are $U$-statistics with symmetric kernels
\begin{align*}
f((x,y),(x^\prime,y^\prime))=\frac{1}{2} \left( 1_{\{x\leq x^\prime   \}}  + 1_{\{x\geq x^\prime   \}}   \right), \
&g((x,y),(x^\prime,y^\prime))=\frac{1}{2} \left( 1_{\{y \leq y^\prime   \}}  + 1_{\{y\geq y^\prime   \}}   \right),\\
h((x,y),(x^\prime,y^\prime))=&\frac{1}{2} \left( 1_{\{x\leq x^\prime, y \leq y^\prime   \}}  + 1_{\{x\geq x^\prime, y\geq y^\prime   \}} \right).
\end{align*}
Note that the underlying random vectors $(\bX_i^{(h)},\bY_i^{(h)})$, $i\geq 1$, are dependent, so that standard $U$-statistics theory for independent data does not apply. However, we can apply an ergodic theorem for  $U$-statistics established in \cite{aaronson:1996}.

\begin{theorem}
Assume that $(X_i,Y_i)$, $i\geq 1$,  is a stationary ergodic process, and that $(\bX^{(h)},\bY^{(h)})$ has  a continuous distribution.
Then, as $n\rightarrow \infty$, we obtain almost surely
\begin{eqnarray*}
 \hat{p}_{\bX^{(h)},n}\longrightarrow p_{\bX^{(h)}},\
  \hat{p}_{\bY^{(h)},n}\longrightarrow p_{\bY^{(h)}},\
 \hat{p}_{(\bX^{(h)},\bY^{(h)}),n}\longrightarrow p_{(\bX^{(h)},\bY^{(h)})}.
\end{eqnarray*}
\end{theorem}
\begin{proof}
We apply Theorem U from \cite{aaronson:1996}. %The conditions of that theorem are satisfied. 
%If the distribution of $(\bX^{(h)},\bY^{(h)})$ is discrete, condition (i) of Theorem U holds. If the distribution is continuous, 
The kernels $f,g$, and $h$ are almost everywhere continuous and thus condition (ii) of Theorem U holds.
\end{proof}

In order to establish asymptotic normality of these estimators, we have to make some assumptions assuring short-range dependence of the underlying process. We will use the concept of near-epoch dependence in probability introduced in \cite{dehling:2017}. This concept is a variation of the usual $L_2$-near-epoch dependence and does not require any moment assumptions. 

\begin{definition}
(i) Given two sub-$\sigma$-fields $\A,\B\subset \F$, we define the absolute regularity coefficient
\[
  \beta(\A,\B)=\sup\{ \sum_{i,j} |\Prob(A_i\cap B_j)-\Prob(A_i)\Prob(B_j)| \},
\] 
where the supremum is taken over all integers $m,n\geq 1$, all partitions  $A_1,\ldots,A_m \in \A$,  and all partitions $B_1,\ldots,B_n \in \B$ of the sample space $\Omega$.
\\[1mm]
(ii) For a stationary stochastic process $Z_i$, $i\in \Z$, we define the absolute regularity coefficients 
\[
  \beta_k:= \beta(\F_{-\infty}^0,\F_k^\infty),
\]
where $\F_k^l$ denotes the $\sigma$-field generated by the random variables $X_k,\ldots,X_l$. The process $(Z_i)_{i\in \Z}$ is called absolutely regular if $\lim_{k\rightarrow \infty}\beta_k=0$. 
\\[1mm]
(iii) An $\R^d$-valued stochastic process $X_i$, $i\geq 1$, is called near-epoch dependent in probability (in short $P$-NED) on the stationary process $Z_i$, $i\in \Z$, if $(X_i,Z_i)$, $i\geq 1$,  is a stationary process, and if there exists a sequence $(a_k)_{k\geq 0}$ of approximating constants with $\lim_{k\rightarrow \infty}a_k=0$,  a sequence of functions $f_k:\R^{2k+1} \rightarrow \R^d$, and a nonincreasing function $\Phi:(0,\infty)\rightarrow (0,\infty)$ such that
\[
  \Prob(|X_0-f_k(Z_{-k},\ldots, Z_k)| \geq \epsilon) \leq a_k \Phi(\epsilon).
\]
\end{definition}

\begin{proposition}
Let $(\bX_i,\bY_i)$,  $i\geq 1$, be a stationary process that is $P$-NED on an absolutely regular process $Z_k$, $k\in \Z$, and assume that 
\[
 a_k\Phi(k^{-6}) =O(k^{-6(2+\delta)/\delta}) \; \mbox{ and } \sum_{k=1}^\infty k\beta_k^{\delta/(2+\delta)}<\infty
\]
for some $\delta>0$. Moreover, assume that  $\bY^{(h)}-\bX^{(h)}$ has a bounded density. Then, the following approximations hold:
\begin{eqnarray*}
\sqrt{n} \left(\hat{p}_{\bX^{(h)},n}- p_{\bX^{(h)}} \right) &=& \frac{2}{\sqrt{n}} 
\sum_{i=1}^n f_1(\bX_i^{(h)},\bY_i^{(h)}  )+o_P(1), \\
\sqrt{n} \left(\hat{p}_{\bY^{(h)},n}- p_{\bY^{(h)}} \right) &=& \frac{2}{\sqrt{n}} 
\sum_{i=1}^n g_1(\bX_i^{(h)},\bY_i^{(h)}  )+o_P(1), \\
\sqrt{n} \left(\hat{p}_{(\bX^{(h)},\bY^{(h)}),n}- p_{(\bX^{(h)},\bY^{(h)})} \right) &=& \frac{2}{\sqrt{n}} 
\sum_{i=1}^n h_1(\bX_i^{(h)},\bY_i^{(h)}  )+o_P(1),
\end{eqnarray*}
where the functions $f_1$, $g_1$, and $h_1$ are the first order terms in the Hoeffding decomposition of the kernels $f$, $g$, and $h$, respectively. 
\label{prop:H-decomp}
\end{proposition}

\begin{remark}
The first order term of the Hoeffding decomposition is given by
\begin{align*}
 f_1(x,y) =& Ef((x,y),(\tilde{\bX}^{(h)},\tilde{\bY}^{(h)}))- p_{\bX^{(h)}}
 = \frac{1}{2} \left( \Prob(\bX^{(h)}\leq x) + \Prob(\bX^{(h)}\geq x)   \right) - \Prob(\bX^{(h)}\leq \tilde{\bX}^{(h)})\\
 =& \frac{1}{2} \left( F(x)+\bar{F}(x) \right) - \Prob(\bX^{(h)}\leq \tilde{\bX}^{(h)}),
\end{align*}
where $F(x):=\Prob(\bX^{(h)} \leq x)$ and $\bar{F}(x):=\Prob(\bX^{(h)}\geq x)$. Similarly, we get 
\begin{eqnarray*}
g_1(x,y) &=& \frac{1}{2} \left( G(x)+\bar{G}(x) \right) - \Prob(\bY^{(h)}\leq \tilde{\bY}^{(h)}), \ \\
h_1(x,y) &=& \frac{1}{2} \left( H(x,y)+\bar{H}(x,y) \right) - \Prob(\bX^{(h)}\leq \tilde{\bX}^{(h)}, \bY^{(h)}\leq \tilde{\bY}^{(h)}),
\end{eqnarray*}
where $G$, $\bar{G}$, $H$, and $\bar{H}$ are defined analogously to $F$ and $\bar{F}$. 
\end{remark}
\begin{proof}[Proof of Proposition~\ref{prop:H-decomp}] This follows from Lemma~D.6 of \cite{dehling:2017} noting that the variation condition is satisfied because the distribution of $\bY^{(h)}-\bX^{(h)}$ has a bounded density. 
\end{proof}

\begin{theorem}\label{thm:sigma}
 Under the same assumptions as in Proposition~\ref{prop:H-decomp}, we have
\[
\sqrt{n} \left( \begin{array}{c} 
\hat{p}_{\bX^{(h)},n}- p_{\bX^{(h)}} \\
\hat{p}_{\bY^{(h)},n}- p_{\bY^{(h)}} \\
\hat{p}_{(\bX^{(h)},\bY^{(h)}),n}- p_{(\bX^{(h)},\bY^{(h)})} 
\end{array}\right)
\stackrel{\mathcal{D}}{\longrightarrow} N(0,\Sigma),
\]
where $\Sigma\in \mathbb{R}^{3\times 3}$ is the limit covariance matrix whose diagonal and off-diagonal entries are given by
\begin{eqnarray*}
 \sigma_{11}&=& \Var\left( f_1(\bX_1^{(h)},\bY_1^{(h)}) \right) + 2\sum_{i=2}^\infty 
 \Cov\left(f_1(\bX_1^{(h)},\bY_1^{(h)}),  f_1(\bX_i^{(h)},\bY_i^{(h)}) \right), \\
 \sigma_{12}&=&\Cov\left( f_1(\bX_1^{(h)},\bY_1^{(h)}), g_1(\bX_1^{(h)},\bY_1^{(h)}) \right) 
  + \sum_{i=2}^\infty  \Cov\left(f_1(\bX_1^{(h)},\bY_1^{(h)}),  g_1(\bX_i^{(h)},\bY_i^{(h)}) \right)\\
 &&  + \sum_{i=2}^\infty  \Cov\left(f_1(\bX_i^{(h)},\bY_i^{(h)}),  g_1(\bX_1^{(h)},\bY_1^{(h)}) \right). 
\end{eqnarray*}

\end{theorem}
\begin{proof}
By the multivariate central limit theorem for partial sums of NED processes we obtain
\[
 S_n:= \frac{1}{\sqrt{n}} \sum_{i=1}^n \left(  f_1(\bX_i^{(h)},\bY_i^{(h)}  ),  g_1(\bX_i^{(h)},\bY_i^{(h)}  ), h_1(\bX_i^{(h)},\bY_i^{(h)}  )  \right)  \stackrel{\mathcal{D}}{\longrightarrow} N(0,\Sigma);
\]
see, e.g.,\cite{wooldridge:1988}.  Now, the statement of the theorem follows from Proposition~\ref{prop:H-decomp} together with an application of Slutsky's lemma.
\end{proof}

\begin{theorem}
Under the  assumptions of Proposition~\ref{prop:H-decomp}, the estimator $\hat{\tau}_n(\bX^{(h)},\bY^{(h)})$ of Kendall's $\tau$ $\tau(\bX^{(h)},\bY^{(h)})$ is consistent and asymptotically normal. More precisely, we obtain
\[
  \sqrt{n} \left( \hat{\tau}_n(\bX^{(h)},\bY^{(h)}) - \tau(\bX^{(h)},\bY^{(h)})  \right) \stackrel{\mathcal{D}}{\longrightarrow} N(0, (\nabla \psi) \Sigma (\nabla \psi)^\top),
\]
where $\psi$ is defined by \eqref{eq:psi} and $\Sigma$ is defined as in Theorem \ref{thm:sigma}.
\end{theorem}
\begin{proof}
This follows from the previous theorem, together with the delta method applied to the function $\psi$.
\end{proof}

\section{Ordinal Pattern Dependence in Contrast to Other Dependence Measures}
\label{sec:comparison}
For independent vectors $(\bX_i,\bY_i)$, $1\leq i\leq n$,
all dependence measures considered in the previous sections make sense. Yet, for measuring dependence between two time series,  only ordinal pattern dependence and Kendall's $\tau$ seem to be reasonable choices of dependence measures. In this section, we 
point out
what kind of dependencies are  measured by ordinal pattern dependence and how  ordinal pattern dependence compares to classical dependence measures such as Pearson's correlation coefficient  and multivariate Kendall's $\tau$. 

\subsection{The case $h=1$}

Axiom (4) in Definition \ref{def:dep_measure} ensures that a multivariate dependence measure takes the value zero if the respective vectors are independent.
In this regard, 
a natural question that arises when studying the dependence between two random vectors is whether the considered dependence measure may also differentiate between independent vectors and uncorrelated, but dependent, random vectors. 
In this section, we provide an answer to this question by giving   examples of marginally uncorrelated Gaussian random vectors with non-vanishing ordinal pattern dependence.
For this purpose, we initially characterize ordinal pattern dependence of order $1$ for Gaussian random vectors. 

\begin{proposition} \label{prop:OPD_1}
Let $\bX=(X_1,X_2)$ and $\bY=(Y_1,Y_2)$ be two Gaussian random vectors satisfying $\E(X_1)=\E(X_2)$, $\E(Y_1)=\E(Y_2)$ and $\Var(X_2-X_1)\neq 0 \neq \Var(Y_2-Y_1)$. Then, it holds that
\begin{align}
\OPD_1(\bX,\bY)  &=\tau(X_2-X_1,Y_2-Y_1)=\frac{2}{\pi} \arcsin \Cor(X_2-X_1, Y_2-Y_1).
\end{align}
\end{proposition}
\begin{proof}
By definition 
\[
\OPD_1(\bX,\bY)=\frac{\Prob(\Pi(X_1,X_2)=\Pi(Y_1,Y_2)) - \sum_{\pi\in S_1}  \Prob(\Pi(X_1,X_2)=\pi)\Prob(\Pi(Y_1,Y_2)=\pi) }{1-\sum_{\pi\in S_1}\Prob(\Pi(X_1,X_2)=\pi)\Prob(\Pi(Y_1,Y_2)=\pi)}.
\]
Since $X_2-X_1$ and $Y_2-Y_1$ are both Gaussian random variables with mean zero and non-zero variance, we obtain
\[
  \Prob(X_1<X_2)=\Prob(X_2<X_1)=\Prob(Y_1<Y_2)=\Prob(Y_2<Y_1)=\frac{1}{2},
\]
and thus $\Prob(\Pi(X_1,X_2)=\pi)=\Prob(\Pi(X_1,X_2)=\pi)=\frac{1}{2}$ for any $\pi\in S_1$. Hence, we obtain
\[
 \OPD_1(\bX,\bY)=2\, \Prob(\Pi(X_1,X_2)=\Pi(Y_1,Y_2)) -1.
\]
Moreover, it holds that
\begin{align*}
\Prob(X_1<X_2, Y_1<Y_2) &=\Prob(X_2-X_1>0, Y_2-Y_1>0)  
=\Prob(X_1-X_2>0, Y_1-Y_2>0)=\Prob(X_1>X_2,Y_1>Y_2),
\end{align*}
and hence 
  \[
   \OPD_1(\bX,\bY)=4\, \Prob( X_1-X_2>0, Y_1-Y_2>0 ) -1.
\]
From Lemma~\ref{lemma1} with $h=1$, we find $\tau(X_2-X_1,Y_2-Y_1)=4 \, \Prob(X_2-X_1<0, Y_2-Y_1<0)-1$, and thus
$\OPD_1(\bX,\bY)=\tau(X_2-X_1,Y_2-Y_1)$. 

Finally, using the  {\em orthant probabilities} formula for Gaussian random variables we obtain
\[
 \Prob(X_2-X_1<0, Y_2-Y_1<0) =\frac{1}{4}+\frac{1}{2\pi} \arcsin \Cor(X_2-X_1,Y_2-Y_1),
\]
and thus $\OPD_1(\bX,\bY) = \frac{2}{\pi}\arcsin \Cor(X_2-X_1,Y_2-Y_1)$.
\end{proof}

In the following we provide an example of marginally uncorrelated Gaussian random vectors
with non-vanishing ordinal pattern dependence of order 1. 

\begin{definition}
A stationary bivariate Gaussian process $W_i=(X_i,Y_i)$ is called $\arp(1)$ process, if there exists a matrix $A\in \R^{2\times 2 } $, and an i.i.d. $N(0,I_2)$-distributed Gaussian process $\xi_i=(\epsilon_i,\eta_i)$  such that the $\arp(1)$-equation 
\begin{equation}
W_{i}=AW_{i-1} +\xi_{i}
\label{eq:biv-ar1}
\end{equation}
is satisfied.
\end{definition}

\begin{remark}
Given a matrix $A\in \R^{2\times 2}$, a stationary Gaussian $\arp(1)$-process exists, if and only if all eigenvalues of $A$ are strictly less than $1$ in absolute value. 
\end{remark}

We will now consider the special example when the $\arp(1)$-matrix is given by
\begin{equation}
 A=\left( \begin{array}{cc}  
  a & b \\
  b & -a
 \end{array}  \right),
 \label{eq:ar1-matrix}
\end{equation}
where $a^2+b^2 <1$. 
Thus, the $\arp(1)$-equation for $W_i=(X_i,Y_i)$ takes the form
\begin{align*}
  X_i =  aX_{i-1} +b Y_{i-1} +\epsilon_i, \
  Y_i = b X_{i-1} -a Y_{i-1} +\eta_i.
\end{align*}
In the next lemma, we state some properties of the processes $(X_i)_{i\geq 1}$ and $(Y_i)_{i\geq 1}$, and we give an explicit formula for their ordinal pattern dependence of order $1$. 

\begin{lemma}
Consider the stationary bivariate Gaussian $\arp(1)$-process $W_i=(X_i,Y_i), i\geq 1$, satisfying  \eqref{eq:biv-ar1} with matrix $A$ given by \eqref{eq:ar1-matrix} such that $a^2+b^2<1$. Then, it holds that
\begin{align}
 \Cov(X_1,Y_1)=0, \
 &\Var(X_1)=\Var(Y_1)=\frac{1}{1-a^2-b^2}  \nonumber \\
 \OPD_1 (\bX^{(1)},\bY^{(1)}) & = \frac{2}{\pi} \arcsin \left(  - \frac{b}{\sqrt{1-a^2}} \right) 
 \label{eq:opd1-ar1}
\end{align}
\label{lem:opd1-ar1}
\end{lemma}
\begin{proof}
The eigenvalues of $A$ are $\lambda_{1,2}=\pm \sqrt{a^2+b^2}$, and thus \eqref{eq:biv-ar1} has a unique stationary solution. Since the $\arp(1)$-equation defines a Markov chain with state space $\R^2$, the joint distribution of $W_i=(X_i,Y_i)$ is uniquely characterized by the distributional fixed point equation
\[
  W\stackrel{\mathcal{D}}{=} A\, W + \xi,
\]
where $\xi =(\epsilon,\eta)$ has a bivariate normal distribution with mean zero and covariance matrix $I_2$, and where $\xi$ is independent of $W$. We will now show that  $W\sim N(0,\sigma^2 I_2)$ satisfies this equation with 
\[
  \sigma^2 = \frac{1}{1-a^2-b^2}. 
\]
In order to prove this, we need to calculate the distribution of $A\, W+\xi$. Since the distribution is Gaussian, it suffices to calculate the variances and the covariance. We obtain
\begin{align*}
&\Cov(aX+bY+\epsilon, bX-aY+\eta)=ab \sigma^2 -ab \sigma^2=0,\\
&\Var(aX+bY+\epsilon)= a^2 \sigma ^2 + b^2 \sigma^2 +1 =\frac{a^2+b^2}{1-a^2-b^2}+1=\frac{1}{1-a^2-b^2}=\sigma^2,\\
& \Var(bX-aY+\epsilon)= b^2 \sigma ^2 + a^2 \sigma^2 +1 =\frac{a^2+b^2}{1-a^2-b^2}+1=\frac{1}{1-a^2-b^2}=\sigma^2,
\end{align*}
which shows that $A\, W+\xi$ has indeed the same distribution as $W$. 
\\[1mm]
In order to determine the $\OPD_1$ of the two processes, we need to calculcate the  correlation of the differences. The covariance of the increments is given by 
\begin{align*}
  \Cov(X_2-X_1, Y_2-Y_1) = \Cov((a-1)X_1+bY_1+\epsilon_2,bX_1 +(-a-1) Y_1 +\eta_2)
  &= b(a-1)\sigma^2 -b(a+1)\sigma^2 \\
  &=-2b\sigma^2 =\frac{-2b}{1-a^2-b^2}
 \end{align*}
 and the variances of the increments are given by
 \begin{align*}
  \Var(X_2-X_1)&= (a-1)^2 \sigma^2 +b^2\sigma^2 +1
     = \frac{ (a-1)^2+b^2 }{1-a^2-b^2} +1 =\frac{2(1-a) }{ 1-a^2-b^2 },\\[2mm]
   \Var(Y_2-Y_1)&= b^2 \sigma^2 +(a+1)^2\sigma^2 +1
   =\frac{b^2+(a+1)^2 }{1-a^2-b^2}+1 = \frac{2(a+1)}{1-a^2-b^2}.
\end{align*}
Thus, we obtain the following formula for the correlation of the increments:
\[
  \Cor(X_2-X_1,Y_2-Y_1) =\frac{-2b}{\sqrt{4(1-a)(a+1)  } }=-\frac{b}{\sqrt{1-a^2}  }
\]
Using the identity $\OPD_1((X_1,X_2),(Y_1,Y_2))=\frac{2}{\pi}\arcsin \Cor(X_2-X_1,Y_2-Y_1)$, we finally obtain \eqref{eq:opd1-ar1}.
\end{proof}

\begin{figure}[h]
\begin{center}
\includegraphics[width=0.3\textwidth]{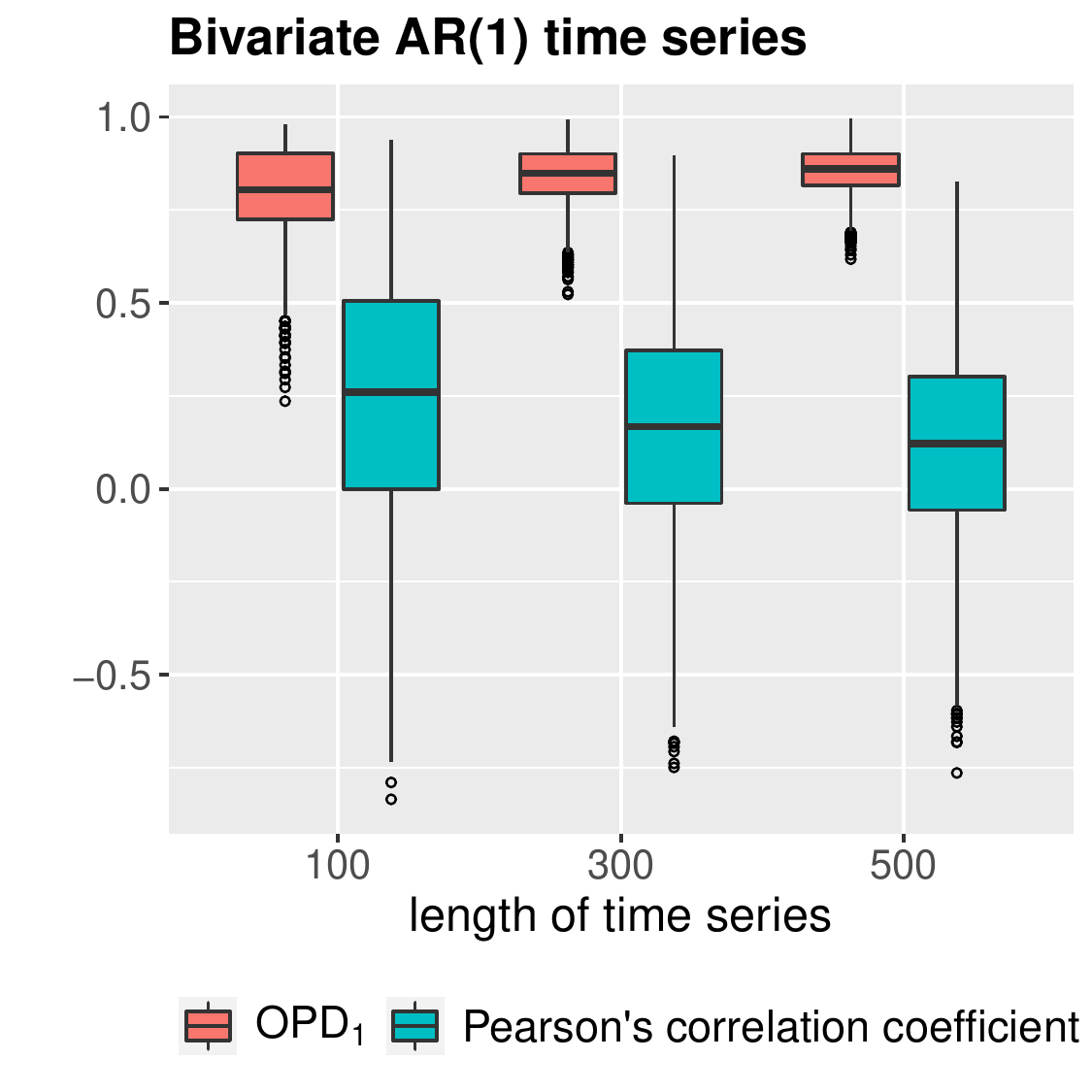}
\caption{
Boxplots of Pearson's correlation coefficient and ordinal pattern dependence of order $h=1$  based on $5000$ repetitions of a bivariate AR(1)-process  $\pr{X_i, Y_i}$, $i\in\{1, \ldots, n\}$, satisfying \eqref{eq:matrix_A} with $a=0.7$ and $b=-0.7$.}
\label{fig:opd_and_pearson}
\end{center}
\end{figure}

\begin{remark}\label{rem:ar1}
(i) The special choice of the $\arp(1)$-matrix  $A$ made in \eqref{eq:ar1-matrix} assures that the two processes $(X_i)_{i\geq 1}$ and $(Y_i)_{i\geq 1}$ have identical marginals, and that $X_i$ and $Y_i$ are independent for each fixed $i$. In fact, one can show that the latter two properties only hold if $A$ is either of the form \eqref{eq:ar1-matrix} or of the form
\begin{align}
  A=\left( \begin{array}{cc}
    a & b\\
    -b & a
  \end{array}  \right). \label{eq:ar1-matrix-new}
\end{align}
In this case, using similar calculations as above, one obtains $\OPD_1((X_1,X_2),(Y_1,Y_2))=0$. 
\\[1mm]
(ii) Lemma~\ref{lem:opd1-ar1} provides an example of a Gaussian process for which Pearson's correlation of $X_i$ and $Y_i$ equals $0$, i.e., the one-dimensional marginals are independent. However, the processes $(X_i)_{\i\geq 1}$ and $(Y_i)_{i\geq 1}$ are not independent, as can be seen from the identity for  $\OPD(\bX^{(1)},\bY^{(1)})$.
\end{remark}
We illustrate our results by simulating a bivariate AR$(1)$-process 
\begin{align*}
W_i\defeq \begin{pmatrix} X_i\\
Y_i\end{pmatrix}, \ i\in \{1, \ldots, 500\},
\end{align*}
with $W_i=AW_{i-1}+\xi_i$, where
\begin{align}\label{eq:matrix_A}
A\defeq\begin{pmatrix} a & b\\
b & -a
\end{pmatrix}, \  \xi_i\defeq \begin{pmatrix} \varepsilon_i\\
\eta_i\end{pmatrix},
\end{align}
$\xi_i$ being a multivariate Gaussian random vector with covariance matrix $\Sigma_{\xi}=I_2$ (with $I_2$ denoting the identity matrix). 
We choose $a^2+b^2<1$, but close to $1$, in order to obtain $\Cov(X_i,Y_i)=0$, but high ordinal pattern dependence.
For the  simulations  summarized by the boxplots in Fig. \ref{fig:opd_and_pearson} we chose $a=0.7$ and $b=-0.7$.
Clearly, the
median of the boxplots that are based on the values of Pearson's correlation coefficient approaches zero, while the median of the boxplots that are based on the values of ordinal pattern dependence seem to converge to a value between $0.75$ and $1$.

\begin{figure}[h]
\begin{center}
\includegraphics[width=0.6\textwidth]{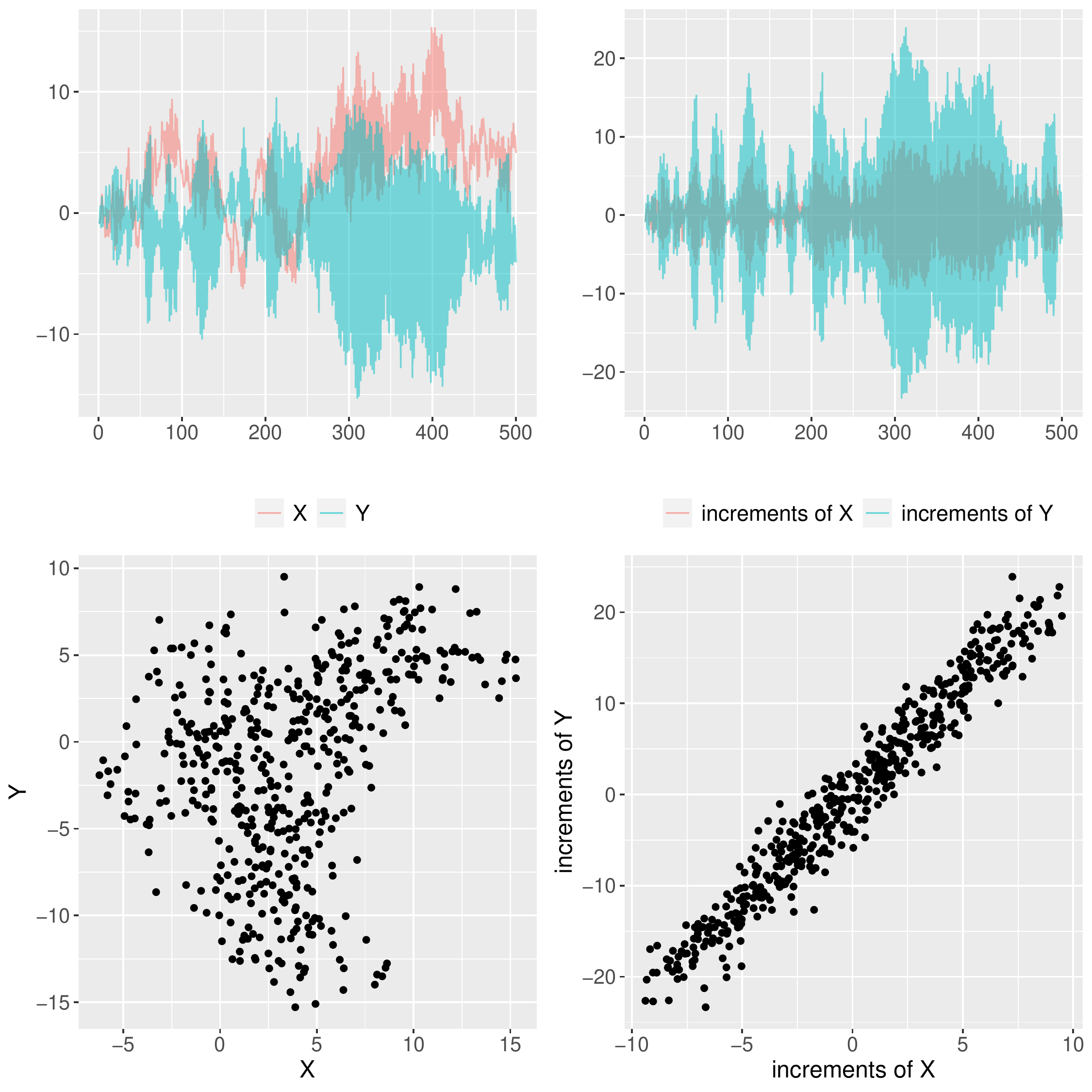}
\caption{
Sample paths  of $X_i$, $i\in\{1, \ldots, 500\}$, and $Y_i$, $i\in\{1, \ldots, 500\}$, their increments at lag $1$,  as well as corresponding scatterplots 
based on a  bivariate AR(1)-process  $\pr{X_i, Y_i}$, $i\in\{1, \ldots, n\}$, satisfying \eqref{eq:matrix_A} with $a=0.7$ and $b=-0.7$.}
\label{fig:scatterplot}
\end{center}
\end{figure}

Fig. \ref{fig:scatterplot} depicts one sample path of the single time series $X_i$, $i\in\{1, \ldots, 500\}$, and $Y_i$, $i\in \{1, \ldots, 500\}$, the corresponding  increment processes, as well as scatterplots of the original observations and their increments. The scatterplots clearly indicate that, while the original observations are uncorrelated,  the increment processes are positively correlated.
Moreover, the scatterplots of the two processes  and their increments in Fig. \ref{fig:scatterplot} underline uncorrelatedness of the original processes and a high dependence of their increments.

\subsection{The case $h=2$}
Let us recall that for the computation of the ordinal pattern dependence of order $h=1$,  
the crucial quantity  is $\Cor(X_2-X_1,Y_2-Y_1)$ since, according to Proposition \ref{prop:OPD_1}, $\OPD_1(\bX^{(1)}, \bY^{(1)})$ is just a monotone transformation of this correlation. It is, therefore, natural to wonder whether it is possible to 
construct a stationary, bivariate process $(X_i,Y_i)_{i\geq 1}$ with $\OPD_1(\bX^{(1)}, \bY^{(1)})=0$, but $\OPD_2(\mathbf{X}^{(2)}, \mathbf{Y}^{(2)})\neq 0$.

The  AR(1)-process in Lemma \ref{lem:opd1-ar1} does not fulfill these conditions, since the restriction 
\begin{align*}
\Cor(X_2-X_1,Y_2-Y_1)
=-\frac{b}{\sqrt{1-a^2}}=0
\end{align*}
implies  $b=0$. As a result, we obtain a process $W_i=\left(X_i,Y_i\right)=\left(aX_{i-1}+\xi_i,-aY_{i-1}+\eta_i\right)$, that does not incorporate any dynamical dependence between the processes $X_i$, $i\geq 1$, and $Y_i$, $i\geq 1$. The only dependence in this model exists within each component. Yet, this  does not have an impact on ordinal pattern dependence.

Following Remark \ref{rem:ar1}, the choice of the matrix $A$ in \eqref{eq:ar1-matrix-new}  yields $\Cor\left(X_i,Y_i\right)=0$ for $i\in\{1, 2\}$, and $\OPD_1(\mathbf{X}^{(1)}, \mathbf{Y}^{(1)})=0$. This leads to the question whether this special construction of an AR(1)-process fulfills $\OPD_2(\mathbf{X}^{(2)}, \mathbf{Y}^{(2)})\neq 0$.
\begin{lemma}
Consider the stationary bivariate Gaussian $\arp(1)$-process $W_i=(X_i,Y_i), i\geq 1$, satisfying  \eqref{eq:biv-ar1} with matrix $A$ given by \eqref{eq:ar1-matrix-new}, and where $a^2+b^2<1$. Then, it holds that
\begin{align*}
 \Cov(X_i,Y_i)=0, \
 &\Var(X_i)=\Var(Y_i)=\sigma^2=\frac{1}{1-a^2-b^2}, \\
 \OPD_1 ((X_i,X_{i+1}), (Y_i,Y_{i+1}))  = 0,\
 \Cor&\left(X_2-X_1,Y_3-Y_2\right)=-b,\
  \Cor\left(X_3-X_2,Y_2-Y_1\right)=b.
\end{align*}
\label{lem:opd2-ar1}
\end{lemma}
\begin{proof}
The first three identities can be shown as in Lemma \ref{lem:opd1-ar1}. Thus, it remains to show the latter two. It holds that
\begin{align*}
\Var(X_2-X_1)&=(a-1)^2\sigma^2+b^2\sigma^2+1
=\frac{(a-1)^2+b^2+1-a^2-b^2}{1-a^2-b^2}
=2(1-a)\sigma^2.
\end{align*}
Analogously, we obtain
\begin{align*}
\Var(Y_3-Y_2)=2(1-a)\sigma^2.
\end{align*}
Furthermore, it holds that
\begin{align*}
\Cov(Y_3-Y_2,X_2-X_1)&=\mathbb{E}\left(Y_3X_2\right)-\mathbb{E}\left(Y_2X_2\right)-\mathbb{E}\left(Y_3X_1\right)+\mathbb{E}\left(Y_2X_1\right)
=2b(a-1)\sigma^2,
\end{align*}
since
\begin{align*}
\mathbb{E}(Y_3X_2)=-b\sigma^2, \
\mathbb{E}(Y_3X_1)=-2ab\sigma^2.
\end{align*}
Alltogether, we arrive at
\begin{align*}
\Cor(Y_3-Y_2,X_2-X_1)=\frac{2ab-2b}{2(1-a)}=-b. 
\end{align*}
$\Cor(Y_2-Y_1,X_3-X_2)=b$ is derived by similar calculations.
\end{proof}
%In order to compute $\OPD_2(\mathbf{X}^{(2)}, \mathbf{Y}^{(2)})$, we have to calculate 
%\begin{align*}
%\Prob\left(\Pi(X_1, X_2, X_{3})=\pi, \Pi(Y_1, Y_2, Y_{3})=\pi\right)
%\end{align*} 
%for every $\pi\in S_{2}$. 
%With $\pi=(0, 1, 2)$ it follows that
%\begin{align*}
%\Prob\left(\Pi(X_1, X_2, X_{3})=\pi, \Pi(Y_1, Y_2, Y_{3})=\pi\right)
%=&\Prob\left(X_1\leq X_2 \leq X_3, Y_1\leq Y_2 \leq Y_3\right)\\
%=&\Prob\left(X_2-X_1 \geq 0, X_3-X_2 \geq 0, Y_2-Y_1\geq 0, Y_3-Y_2\geq 0\right).
%\end{align*}
%As a result, computing $\OPD_2(\mathbf{X}^{(2)}, \mathbf{Y}^{(2)})$ boils down to determining the orthant probabilities of a four-dimensional Gaussian vector. To our knowledge  a closed expression for these probabilities
%is not at hand; see \cite{abrahamson:1964}. However, the methods introduced in  \cite{abrahamson:1964} show that the correlation of the increments at lag $1$ is a substantial parameter to determine the orthant probabilities - and, therefore, $\OPD_2$. Before emphasizing this empirically by simulations, we construct AR(h)-processes to obtain $\OPD_h\left( \bX^{(h)},\bY^{(h)} \right)\neq 0$ but $\OPD_i\left( \bX^{(i)},\bY^{(i)} \right)=0$, $i=1,\ldots,h-1$ and $\mathrm{Cor}\left(X_1,Y_1\right)=0$. Note that $\OPD_h\left( \bX^{(h)},\bY^{(h)} \right)$ is closely connected to $\Cor(Y_{h+1}-Y_h,X_2-X_1)$ and $\Cor(X_{h+1}-X_h,Y_2-Y_1)$. We focus on $h=2$:

Lemma \ref{lem:opd2-ar1} provides an example of a bivariate process $(X_i,Y_i)_{i\geq 1}$ for which $\mathrm{Corr}(X_i,Y_i)=0$ and $\OPD_1(\bX^{(1)},\bY^{(1)})=0$, but where the processes $(X_i)_{i\geq 1}$ and $(Y_i)_{i\geq 1}$ are nevertheless dependent. The fact that the increments $X_2-X_1$ and $Y_3-Y_2$ are dependent, leads us to conjecture that $\OPD_2(\bX^{(2)},\bY^{(2)})\neq 0$, but we do not have a proof. In order to compute $\OPD_2(\bX^{(2)},\bY^{(2)})$, we have to calculate 
\begin{align*}
\Prob\left(\Pi(X_1, X_2, X_{3})=\pi, \Pi(Y_1, Y_2, Y_{3})=\pi\right)
\end{align*} 
for any $\pi\in \mathcal{S}_2$. This requires computations of orthant probabilities for 4-dimensional Gaussian vectors, e.g. for $\pi=(0,1,2)$ we obtain 
\begin{align*}
\Prob\left(\Pi(X_1, X_2, X_{3})=\pi, \Pi(Y_1, Y_2, Y_{3})=\pi\right)
=&\Prob\left(X_1\leq X_2 \leq X_3, Y_1\leq Y_2 \leq Y_3\right)\\
=&\Prob\left(X_2-X_1 \geq 0, X_3-X_2 \geq 0, Y_2-Y_1\geq 0, Y_3-Y_2\geq 0\right).
\end{align*}
To the best of our knowledge, there are no explicit  formulas for these probabilities known. 

In what follows, we present an example of a bivariate AR(2)-process  $(X_i,Y_i)_{i\geq 1}$ for which $\mathrm{Corr}(X_i,Y_i)=0$ and $\OPD_1(\bX^{(1)},\bY^{(1)})=0$, but where the processes $(X_i)_{i\geq 1}$ and $(Y_i)_{i\geq 1}$ are dependent. For this example, we show by means of a Monte Carlo simulation, that $\OPD_2(\bX^{(2)},\bY^{(2)})\neq 0$.

\begin{example}
Let $W_i$, $i\geq 1$, be a bivariate AR(2)-process defined by
\begin{align*}
W_i\defeq \begin{pmatrix} X_i\\
Y_i\end{pmatrix}, \ i\geq 1,
\end{align*}
where $W_i=AW_{i-2}+\xi_i$ with
\begin{align}\label{eq:matrix_A_2}
A\defeq\begin{pmatrix} a & b\\
b & -a
\end{pmatrix} \ \text{and} \ \xi_i\defeq \begin{pmatrix} \varepsilon_i\\
\eta_i\end{pmatrix},
\end{align}
$\xi_i$,  $i\geq 1$,  bivariate Gaussian random vectors with covariance matrix $\Sigma_{\xi}=I_2$ (with $I_2$ denoting the identity matrix) and
$X_1\defeq \xi_1$, $Y_1=\eta_1$, $X_2\defeq \xi_2$, $Y_2=\eta_2$.
Moreover, we assume that  
   $\sigma^2\defeq\Var(X_1)=\Var(Y_1)$.
By definition it holds that $\Cov(X_1,Y_1)=\Cov(X_2,Y_2)=0$.
Moreover, we have	 $\OPD_1(\mathbf{X}^{(1)}, \mathbf{Y}^{(1)})=0$ since
\begin{align*}
\Cov(X_3-X_2,Y_3-Y_2)&=\mathbb{E}\left[\left(aX_1+bY_1+\xi_3-X_2\right)\left(bX_1-aY_1+\eta_3-Y_2\right)\right]=ab\sigma^2-ba\sigma^2=0.
\end{align*}
This construction of AR(2)-processes can be extended to AR(h) for $h\in\mathbb{N}$, if one wants to obtain $\OPD_h\left( \bX^{(h)},\bY^{(h)} \right)\neq 0$ but $\OPD_i\left( \bX^{(i)},\bY^{(i)} \right)=0$, $i=1,\ldots,h-1$ and $\mathrm{Corr}\left(X_1,Y_1\right)=0$ by using $h$ independent AR(1)-processes and couple them via 
\begin{align*}
X_j=A\left(X^{(1)}_{j-h},X^{(2)}_{j-h}\right).
\end{align*}
To illustrate the strong connection between a large correlation of the increments at lag $1$ and $\OPD_2(\mathbf{X}^{(2)}, \mathbf{Y}^{(2)})$ we consider estimated values for $\OPD_2(\mathbf{X}^{(2)}, \mathbf{Y}^{(2)})$ based on simulations of a bivariate AR$(2)$-process
that satisfies the above assumptions; see Fig. \ref{fig:opd1_and_opd2}.
The corresponding  boxplots clearly indicate that,  as the length of the time series increases, the ordinal pattern dependence of order $h=1$ approaches zero, while the ordinal pattern dependence of order $h=2$ converges to a a value between $0.1$ and $0.25$.

\begin{figure}[h]
\begin{center}
\includegraphics[width=0.3\textwidth]{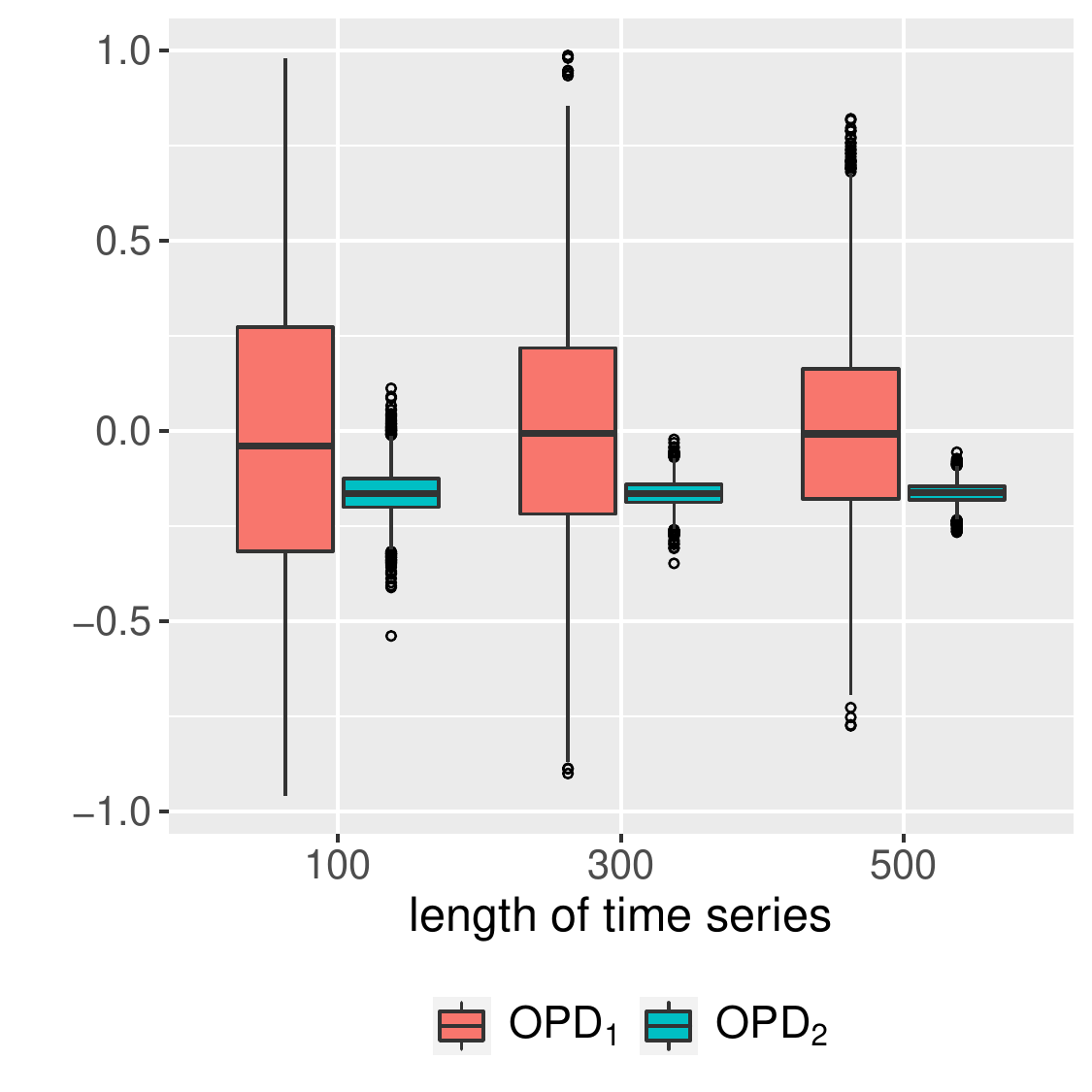}
\caption{Boxplots of ordinal pattern dependence of order $h=1$ and $h=2$ based on $5000$ repetitions of a bivariate AR(2)-process  $\pr{X_i, Y_i}$, $i\in\{1, \ldots, n\}$, satisfying  \eqref{eq:matrix_A_2} with $a=0.01$ and $b=0.98$.}
\label{fig:opd1_and_opd2}
\end{center}
\end{figure}
\end{example}

\subsection{Ordinal pattern dependence in contrast to multivariate Kendall's $\tau$}

Recall that for Gaussian random vectors $\bX^{(1)}=(X_1, X_2)$ and $\bY^{(1)}=(Y_1, Y_2)$ satisfying $\E(X_1)=\E(X_2)$, $\E(Y_1)=\E(Y_2)$, and $\Var(X_2-X_1)\neq 0 \neq \Var(Y_2-Y_1)$, it holds that
\begin{align*}
\OPD_1(\bX^{(1)},  \bY^{(1)})= \tau\pr{\tilde{X}_1, \tilde{Y}_1},
\end{align*}
where $\tilde{X}_1\defeq X_2-X_1$ and  $\tilde{Y}_1\defeq Y_2-Y_1$; see Proposition  \ref{prop:OPD_1}.

In general, the following proposition establishes a relation between ordinal pattern dependence and  multivariate Kendall's $\tau$  for Gaussian processes:

\begin{proposition}\label{prop:OPD_and_Kendall}
Let $(X_i,Y_i)$, $i\geq 1$,  be a bivariate, mean zero stationary Gaussian process. Then, it holds that
\begin{align*}
\OPD_h\left( \bX^{(h)} , \bY^{(h)} \right)
=&\frac{2\tau\left(\tilde{X}_{1},\ldots,\tilde{X}_{h},\tilde{Y}_{1},\ldots,\tilde{Y}_{h}\right)\sqrt{ \tilde{p}_{\tilde{\bX}^{(h)}} \left(1-\tilde{p}_{\tilde{\bX}^{(h)}}\right)\tilde{p}_{\tilde{\bY}^{(h)}}\left(1-\tilde{p}_{\tilde{\bY}^{(h)}}\right)}}{1-\sum_{\pi\in S_h}\Prob\left(\Pi\left( \bX^{(h)}\right)=\pi\right)\Prob\left(\Pi\left( \bY^{(h)}\right)=\pi\right)}\\
&+\sum_{\pi\in S_h\setminus\ T_h}\frac{ \Prob\left(\Pi\left( \bX^{(h)}\right)=\Pi\left( \bY^{(h)}\right)=\pi\right)-\Prob\left(\Pi\left( \bX^{(h)}\right)=\pi\right)\Prob\left(\Pi\left( \bY^{(h)}\right)=\pi\right)}{1-\sum_{\pi\in S_h}\Prob\left(\Pi\left( \bX^{(h)}\right)=\pi\right)\Prob\left(\Pi\left( \bY^{(h)}\right)=\pi\right)},
\end{align*}
where $T_h\defeq \{(h,\ldots,0),(0,\ldots,h)\}$, $\tilde{X_i}=X_{i+1}-X_{i}$, $\tilde{Y_i}=Y_{i+1}-Y_{i}$,
\begin{align*}
\tilde{p}_{\tilde{\bX}^{(h)}}&=\Prob\left(\Pi\left(X_1,\ldots,X_{1+h}\right)=(0,\ldots,h)\right)=\Prob\left(\tilde{X}_{1}\leq 0,\ldots,\tilde{X}_{h}\leq 0\right),\\
\tilde{p}_{\tilde{\bY}^{(h)}}&=\Prob\left(\Pi\left(Y_1,\ldots,Y_{1+h}\right)=(0,\ldots,h)\right)=\Prob\left(\tilde{Y}_{1}\leq 0,\ldots,\tilde{Y}_{h}\leq 0\right),
\end{align*}
and
\begin{align*}
\OPD_h\left( \bX^{(h)} , \bY^{(h)}\right)
=&\frac{ \sum_{\pi \in S_h} \tau\left(  X_{1+\pi_1}- X_{1+\pi_0},\ldots, Y_{1+\pi_{h}}-Y_{1+\pi_{h-1}} \right)\sqrt{p_{\bX^{(h)},\pi}\left(1-p_{\bX^{(h)},\pi}\right)p_{\bY^{(h)},\pi}\left(1-p_{\bY^{(h)},\pi}\right)}}{1-\sum_{\pi\in S_h} p_{\bX^{(h)},\pi}p_{\bY^{(h)},\pi}}
\end{align*}
with 
\begin{align*}
p_{\bX^{(h)},\pi}=\Prob\left(\Pi\left(X_1,\ldots,X_{1+h}\right)=\pi\right), \
p_{\bY^{(h)},\pi}=\Prob\left(\Pi\left(Y_1,\ldots,Y_{1+h}\right)=\pi\right).
\end{align*}
\end{proposition}

\begin{remark}
It is a characteristic of Gaussian observations $(X_i)$, $i\geq 1$, and $(Y_i)$, $i\geq 1$, that the distribution of 
\begin{align*}
\left(X_{1+\pi_1}-X_{1+\pi_0},\ldots,X_{1+\pi_{h}}-X_{1+\pi_{h-1}},Y_{1+\pi_1}-Y_{1+\pi_0},\ldots,Y_{1+\pi_{h}}-Y_{1+\pi_{h-1}}\right) 
\end{align*}
is uniquely determined by the autocovariances and crosscovariances of $\bX^{(h)}$ and $\bY^{(h)}$. For this reason, it is possible to express all  the dependencies in the vector above by the two-dimensional marginals of a multivariate Gaussian distribution. However, since we do not have a closed expression for orthant probabilities of a multivariate Gaussian vector with more than $3$ elements, it is not possible to constitute a closed form for ordinal pattern dependence in terms of Kendall's $\tau$.
\end{remark}

 The proof of Proposition \ref{prop:OPD_and_Kendall} has been postponed to Section \ref{sec:proofs}.
In order to illustrate the relation of multivariate Kendall's $\tau$ and ordinal pattern dependence, characterized through Proposition \ref{prop:OPD_and_Kendall},
we consider the case $h=1$ in the following example:
\begin{example}
Let $\bX=(X_1, X_2)$ and $\bY=(Y_1, Y_2)$ be
Gaussian random vectors as in Proposition \ref{prop:OPD_1}. Recall that
\begin{align*}
\OPD_1(\bX, \bY)=\tau\left(X_2-X_1,Y_2-Y_1\right)=\frac{2}{\pi}\arcsin\left(\Cor(X_2-X_1, Y_2-Y_1)\right).
\end{align*}
Moreover, if $\bX$ and $\bY$ have standard normal marginal distributions, it holds that
\begin{align*}
\Cor(X_2-X_1, Y_2-Y_1)=\frac{2\mathbb{E}\left(X_1Y_1\right)-\mathbb{E}\left(X_1Y_2\right)-\mathbb{E}\left(Y_1X_2\right)}{2-2\mathbb{E}\left(X_1X_2\right)}.
\end{align*}
In general, we know that   $\tau(X,Y)=\frac{2}{\pi}\arcsin\left(\Cor(X,Y)\right)$ for a Gaussian random vector $\left(X,Y\right)$, and hence $\Cor(X,Y)=\sin\left(\frac{\pi}{2}\tau(X,Y)\right)$.\newline
As a result, we obtain
\begin{align*}
\OPD_1\left( \bX, \bY\right)
&=\frac{2}{\pi}\arcsin\left( \frac{2\sin\left(\frac{\pi}{2}\tau(X_1,Y_1)\right)-\sin\left(\frac{\pi}{2}\tau(X_1,Y_2)\right)-\sin\left(\frac{\pi}{2}\tau(X_2,Y_1)\right)}{2-2\sin\left(\frac{\pi}{2}\tau(X_1,X_2)\right)}  \right).
\end{align*}
Therefore, the ordinal pattern dependence of order $1$ is determined  by $\tau\left(X_1,Y_1\right)$, $\tau\left(X_1,Y_2\right)$, $\tau\left(X_2,Y_1\right)$, and $\tau\left(X_2,Y_2\right)$.
\end{example}

\subsection{Simulation Study}

\begin{figure}[h!]
\begin{center}
\includegraphics[width=0.9\textwidth]{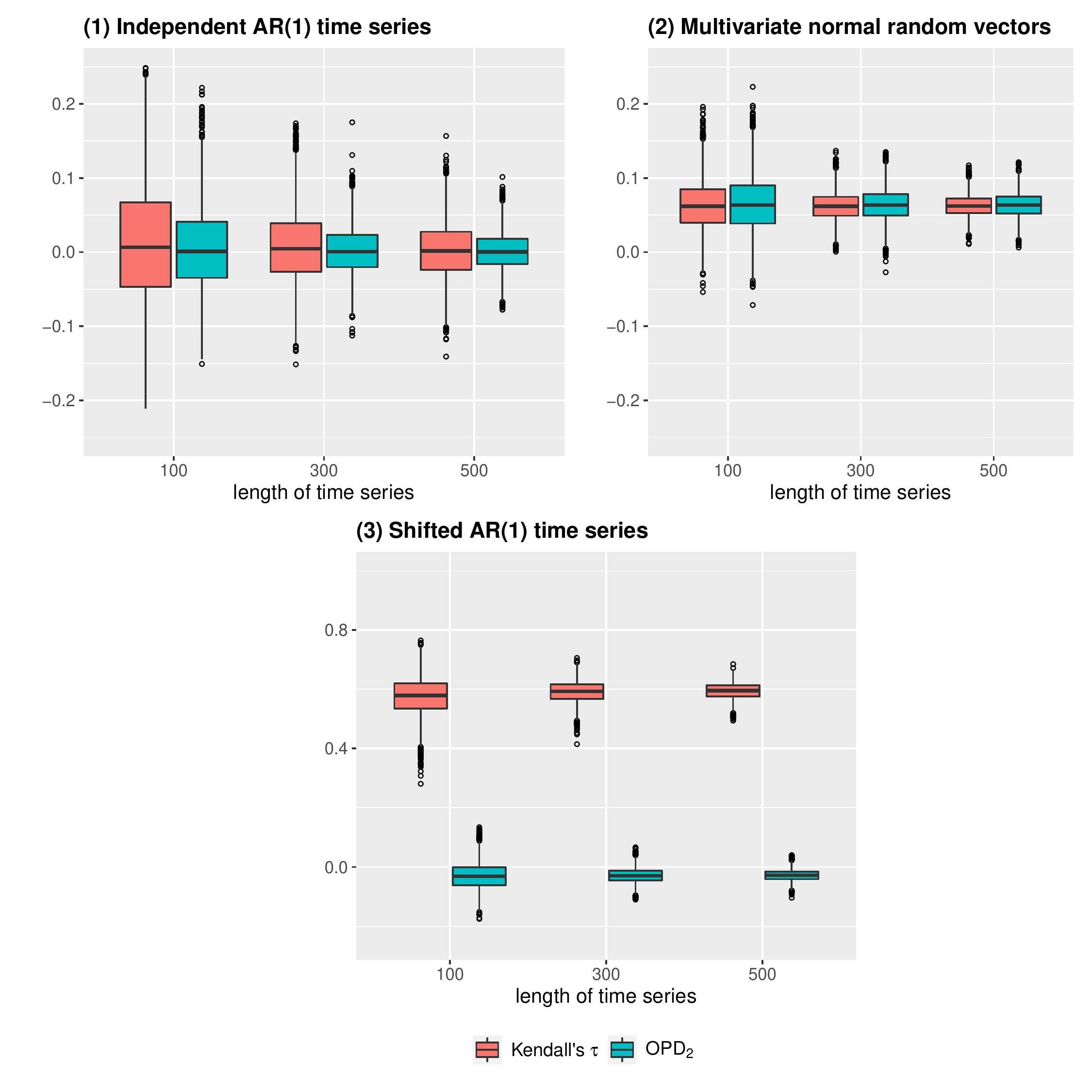}
\caption{Boxplots of ordinal pattern dependence of order  $h=2$ based on $5000$ repetitions of a bivariate process  $\pr{X_i, Y_i}$, $i\in \{1, \ldots, n\}$, corresponding to the situations (1)-(3).}
\label{fig:opd_and_kendall}
\end{center}
\end{figure}

In this section, 
we compare the estimators for ordinal pattern dependence and multivariate Kendall's $\tau$ based on the vectors   $\bX_i^{(h)},\bY_i^{(h)}$
, generated by
 bivariate processes  $(X_i,Y_i)$, $i\geq 1$, in a simulation study. Fig. \ref{fig:opd_and_kendall} corresponds to the case $h=2$. 
 The following situations are considered:
\begin{enumerate}
\item We simulate $(X_i)$, $i\geq 1$, and $(Y_i)$, $i\geq 1$,  as two independent AR(1) time series with
\begin{align*}
X_i=\rho X_{i-1}+\varepsilon_i, \ Y_i=\rho Y_{i-1}+\eta_i,
\end{align*}
where $\left|\rho\right|<1$, and $(\varepsilon_i)$, $i\geq 1$,  and $(\eta_i)$, $i\geq 1$, are two independent  sequences of random variables, both i.i.d. standard normally distributed. In this case, the  sequences $(X_i)$, $i\geq 1$, and $(Y_i)$, $i\geq 1$, are generated by the function
\verb$arima.sim$ in \verb$R$. For the simulations depicted in Fig. \ref{fig:opd_and_kendall}, we chose $\rho=0.5$.
As expected, the values of both dependence measures vary around $0$. Moreover, the boxplots become narrower as the sample size increases confirming consistency of the estimators. The boxplots that correspond to the estimate for Kendall's $\tau$ are wider. This indicates a faster convergence of the estimators for ordinal pattern dependence.
For a systematic simulation study see Table \ref{table:independent} in the appendix.
\item  We simulate $(X_i)$, $i\geq 1$, and $(Y_i)$, $i\geq 1$, as  sequences of independent,  multivariate normal random vectors with values in $\mathbb{R}^3$ and a joint normal distribution with expectation $0$ and covariance matrix
\begin{align}\label{eq:cov}
\Sigma=
\begin{pmatrix}
1 & 0 & 0 & \rho & \rho & \rho  \\
                  0& 1& 0& \rho & \rho & \rho  \\
                  0& 0& 1& \rho & \rho & \rho \\
                  \rho & \rho & \rho & 1& 0& 0 \\
                      \rho & \rho & \rho  & 0& 1& 0 \\
                \rho & \rho & \rho  & 0& 0& 1
\end{pmatrix}.
\end{align}
%$\Sigma$ is positive definite if $\rho\in\{0.1, 0.2, 0,3\}$.
The $\mathbb{R}^6$-valued random vectors  $(X_i,Y_i)$, $i\geq 1$, are generated by the function\lstinline$ rmvnorm$  in \lstinline$R$.
For the simulations depicted in Fig. \ref{fig:opd_and_kendall}, we chose $\rho=0.2$. Note that
the values of both dependence measures deviate from $0$, thereby indicating a correlation between the two processes $(X_i)$, $i\geq 1$, and $(Y_i)$, $i\geq 1$. In fact, the boxplots of both estimators look very similar so that the rates of convergence  seem to be comparable. 
For a systematic simulation study see Table \ref{table:mult_normal} in the appendix.
\item We simulate $(X_i)$, $i\geq 1$, as an AR(1) time series, while $(Y_i)$, $i\geq 1$, corresponds to $(X_i)$, $i\geq 1$, shifted by one time point. More precisely, we simulate 
\begin{align*}
X_i=\rho X_{i-1}+\varepsilon_i, 
\end{align*}
where $\left|\rho\right|<1$, and  $(\varepsilon_i)$, $i\geq 1$, is an i.i.d. standard normally distributed sequence of random variables, and we define
 $Y_i=X_{i+1}$. For the simulations depicted in Fig. \ref{fig:opd_and_kendall}, we chose $\rho=0.5$.

It is intriguing  that ordinal pattern dependence does not detect the high correlation of the time series.  The theoretical value of ordinal pattern dependence  for $h=1$ in this case is given by
\begin{align*}
\OPD_1(\bX^{(1)},\bY^{(1)})&=\frac{2}{\pi} \arcsin \Cor(X_2-X_1,Y_2-Y_1)
 =\frac{2}{\pi} \arcsin \Cor(X_2-X_1,X_3-X_2).
\end{align*}

Routine calculations 
yield the following formula for the ordinal pattern dependence of order 1 between an AR(1) process and the same process shifted by one time point:
\[
\OPD_1(\bX^{(1)},\bY^{(1)}) = \frac{2}{\pi} \arcsin\left(\frac{\rho-1}{2}\right).
\]
As a result,  $\OPD_1(\bX^{(1)},\bY^{(1)})=-0.297, -0.161, -0.032$ 
for   $\rho=0.1, 0.5, 0.9$.
These values coincide with the results of the corresponding  systematic simulation study; see Table \ref{table:shift} in the appendix.

  In \cite{schnurr:dehling:2017} on page 713, an approach is presented to solve the insensitivity of ordinal pattern dependence concerning time shifts. The authors introduced time shifted ordinal pattern dependence in order to investigate and compare time series that are known to have a similar behavior within a certain time deviation. These time series arise for example in the context of hydrology, if discharge data of a river is considered for two different locations. For a real-world data analysis see \cite{nuessgen:2021}. Using this approach, (time-shifted) ordinal pattern dependence of $1$ is detected, since all patterns coincide if we reshift the second time series.
\end{enumerate}

\pagebreak
\section{Proofs}
\label{sec:proofs}

\begin{proof}[Proof of Theorem \ref{thm:opd_dep_measure}]
The first four axioms in Definition \ref{def:dep_measure} are easily verified, the fifth one is involved. We show that the fifth axiom is fulfilled for $\OPD_h$ with $h=2$. For $h=1$,  the difficulties in the proof are not revealed, while for $h>2$ the proof works analogously, but is  notationally more complicated. 
Due to stationarity, it is enough to focus on the first three components of $\mathbf{X}$, $\mathbf{Y}$,  $\mathbf{X^*}$, and $ \mathbf{Y^*}$, i.e., without loss of generality we consider
\[
  \mathbf{X}=(X_1,X_2,X_3), \mathbf{Y}=(Y_1,Y_2,Y_3),  \mathbf{X^*}=(X_1^*,X_2^*,X_3^*), \mathbf{Y^*}=(Y_1^*,Y_2^*,Y_3^*).
\]
Moreover, we can restrict our considerations to $\Prob\left(\Pi\pr{\bX}=\Pi\pr{\bY}\right)$ (the remaining summands  of $\OPD_h(\mathbf{X}	,\mathbf{Y})$ only relate to the distribution of $\mathbf{X}$ and $\mathbf{Y}$ separately). By Axiom 2 (invariance under permutation) it is, furthermore, enough to consider the monotone increasing pattern, that is, $\Pi(X_1, X_2, X_3)=(3, 2, 1)$. Let 
\begin{align}\label{eq:concordance}
\mathbf{X}\overset{\mathcal{D}}{=}\mathbf{X^*}, \ \mathbf{Y}\overset{\mathcal{D}}{=}\mathbf{Y^*}, \ \text{and} \
\binom{\mathbf{X}} {\mathbf{Y}} \preccurlyeq_C  \binom{\mathbf{X^*}} {\mathbf{Y^*}}.
\end{align}
 It is a well-known fact  that \eqref{eq:concordance} implies 
\begin{align*}
\binom{\mathbf{X}} {\mathbf{Y}}^I \preccurlyeq_C  \binom{\mathbf{X^*}}{\mathbf{Y^*}}^I
\end{align*} 
 for all subvectors of variables with indices in $I\subseteq \{1, \ldots, 6\}$, i.e., removing dimensions does not influence which scenario has the larger dependence measure; see \cite{joe:1990}. We will make extensive use of this fact in what follows. 
Moreover, recall that we assume that all considered marginal distribution functions are continuous.
 
 Defining
\[
	\Prob^{x_1,y_1}(A)\defeq \Prob(A|X_1=x_1,Y_1=y_1), 
\]
 considering $\Pi(X_1, X_2, X_3)=(3, 2, 1)$, and using disintegration twice yields
\begin{align*}
&\hspace*{-5mm}\Prob\left(\Pi(X_1, X_2, X_{3})=\Pi(Y_1, Y_2, Y_{3})\right)
=\Prob\pr{X_1\leq X_2\leq X_3, Y_1\leq Y_2\leq Y_3}\\
=&\Prob\pr{\{X_1\leq X_2\}\cap \{X_2\leq X_3\} \cap \{Y_1\leq Y_2\} \cap \{Y_2\leq Y_3\}} 
=\int_{\R^2} \Prob^{x_1,y_1}(x_1\leq X_2 \leq X_3, y_1\leq Y_2\leq Y_3) \ dP_{\binom{X_1}{Y_1}} (x_1,y_1)\\
=&\int_{\R^2} \Prob^{x_1,y_1}(X_2\leq X_3, Y_2\leq Y_3|x_1\leq X_2, y_1\leq Y_2)\cdot \Prob^{x_1,y_1}(X_2\geq x_1, Y_2\geq y_1)    \ dP_{\binom{X_1}{Y_1}} (x_1,y_1) \\
=&\int_{\R^2}  \int_{[x_1,\infty[\times [y_1,\infty[}  \Prob^{x_1,y_1}(X_2\leq X_3, Y_2\leq Y_3| X_2=x_2, Y_2=y_2) \ dP_{\binom{X_2}{Y_2}}^{x_1,y_1}(x_2,y_2)  \Prob^{x_1,y_1}(X_2\geq x_1, Y_2\geq y_1)          \ dP_{\binom{X_1}{Y_1}} (x_1,y_1).
\end{align*}
Since  $\bar{F}_{\binom{X_3}{Y_3}}(x_2,y_2)=\Prob^{x_1,y_1}(X_2\leq X_3, Y_2\leq Y_3|x_1\leq X_2, y_1\leq Y_2)$, it follows that
\begin{align*}
&\Prob\left(\Pi(X_1, X_2, X_{3})=\Pi(Y_1, Y_2, Y_{3})\right)
=&\int_{\R^2}  \int_{[x_1,\infty[\times [y_1,\infty[}   \bar{F}_{\binom{X_3}{Y_3}}(x_2,y_2)  \ dP_{\binom{X_2}{Y_2}}^{x_1,y_1}(x_2,y_2)     \bar{F}_{\binom{X_2}{Y_2}}(x_1,y_1)  \ dP_{\binom{X_1}{Y_1}} (x_1,y_1).
\end{align*}

Due to \eqref{eq:concordance}, we  deduce that
\[
  \int 1_{[a,\infty[\times [b,\infty[}(x,y) \ dP_{\binom{X_j}{Y_j} } (x,y)\leq   \int 1_{[a,\infty[\times [b,\infty[}(x,y) \ dP_{\binom{X^*_j}{Y^*_j} } (x,y)
\]
for $a,b\in\R$ and $j\in\{1,2,3\}$. Since survival functions can be approximated by sums of indicator functions, the bounded convergence theorem yields
\begin{align*}
  \int_{[x_1,\infty[\times [y_1,\infty[}\bar{F}_{\binom{X_3^*}{Y_3^*}}(x_2,y_2)\ dP_{\binom{X_2}{Y_2} }^{x_1,y_1} (x_2,y_2)
  \leq   \int_{[x_1,\infty[\times [y_1,\infty[}\bar{F}_{\binom{X_3^*}{Y_3^*}}(x_2,y_2)\ dP_{\binom{X^*_2}{Y^*_2} }^{x_1,y_1} (x_2,y_2).
\end{align*}
Moreover,  up to scaling, the function 
\[
H(x_1,x_2)\defeq \int_{[x_1,\infty[\times [y_1,\infty[}\bar{F}_{\binom{X_3^*}{Y_3^*}}(x_2,y_2)\ dP_{\binom{X^*_2}{Y^*_2} }^{x_1,y_1} (x_2,y_2)
\]
can be considered as a survival function.
 Hence, we can use an approximation by sums of indicator functions for both,  $H$ and  $\bar{F}_{\binom{X_2^*}{Y_2^*}}$.  
Most notably, 
the product of  two functions having this property is of the same type, i.e., it can as well be approximated by sums of indicator functions. 

Thus, we finally arrive at
\begin{align*}
&\hspace*{-10mm}\Prob\left(\Pi(X_1, X_2, X_{3})=\Pi(Y_1, Y_2, Y_{3})=\left(3, 2, 1\right)\right)\\
\leq &\int_{\R^2} \int_{[x_1,\infty[\times [y_1,\infty[}   \bar{F}_{\binom{X^*_3}{Y^*_3}}(x_2,y_2)  \ dP_{\binom{X^*_2}{Y^*_2}}^{x_1,y_1}(x_2,y_2)       \bar{F}_{\binom{X^*_2}{Y^*_2}}(x_1,y_1)     \ dP_{\binom{X^*_1}{Y^*_1}} (x_1,y_1) \\
=&\Prob(X_1^*\leq X_2^*\leq X_3^*, Y_1^*\leq Y_2^*\leq Y_3^*)=\Prob\left(\Pi(X_1^*, X_2^*, X_{3}^*)=\Pi(Y_1^*, Y_2^*, Y_{3}^*)=(3, 2, 1)\right).
\end{align*}
\end{proof}

\begin{proof}[Proof of Proposition \ref{prop:OPD_and_Kendall}]
First, note that
\begin{align*}
\OPD_h\left( \bX^{(h)} , \bY^{(h)} \right)
=&\frac{\sum_{\pi\in T_h} \left\{\Prob\left(\Pi\left( \bX^{(h)}\right)=\Pi\left( \bY^{(h)}\right)=\pi\right)-\Prob\left(\Pi\left( \bX^{(h)}\right)=\pi\right)\Prob\left(\Pi\left( \bY^{(h)}\right)=\pi\right)\right\}}{1-\sum_{\pi\in S_h}\Prob\left(\Pi\left( \bX^{(h)}\right)=\pi\right)\Prob\left(\Pi\left( \bY^{(h)}\right)=\pi\right)} \\
&+\frac{\sum_{\pi\in S_h\setminus T_h}\left\{ \Prob\left(\Pi\left( \bX^{(h)}\right)=\Pi\left( \bY^{(h)}\right)=\pi\right)-\Prob\left(\Pi\left( \bX^{(h)}\right)=\pi\right)\Prob\left(\Pi\left( \bY^{(h)}\right)=\pi\right)\right\}}{1-\sum_{\pi\in S_h}\Prob\left(\Pi\left( \bX^{(h)}\right)=\pi\right)\Prob\left(\Pi\left( \bY^{(h)}\right)=\pi\right)}.
\end{align*}
Focusing on the pattern $\pi=(0,\ldots,h)$ in the first summand yields
\begin{align*}
\Prob\left(\Pi\left( \bX^{(h)}\right)=\Pi\left( \bY^{(h)}\right)=(0,\ldots,h)\right)
=&\Prob\left(X_1\geq \ldots\geq X_{1+h},Y_1\geq  \ldots\geq Y_{1+h}\right)\\
=&\Prob\left(\tilde{X}_{1}\leq 0,\ldots,\tilde{X}_{h}\leq 0,\tilde{Y}_{1}\leq 0,\ldots,\tilde{Y}_{h}\leq 0  \right)\\
=& \tau\left(\tilde{X}_{1},\ldots,\tilde{X}_{h},\tilde{Y}_{1},\ldots,\tilde{Y}_{h}\right) \sqrt{\tilde{p}_{\tilde{\bX}^{(h)}} \left(1-\tilde{p}_{\tilde{\bX}^{(h)}}\right)\tilde{p}_{\tilde{\bY}^{(h)}}\left(1-\tilde{p}_{\tilde{\bY}^{(h)}}\right)}+\tilde{p}_{\tilde{\bX}^{(h)}}\tilde{p}_{\tilde{\bY}^{(h)}}.
\end{align*}
Due to  symmetry  of the multivariate normal distribution, we have $\left( \bX^{(h)},\bY^{(h)}\right)\overset{\mathcal{D}}{=}\left( -\bX^{(h)},-\bY^{(h)}\right)$.
Therefore, it follows that 
\begin{align*}
\Prob\left(\Pi\left( \bX^{(h)}\right)=\Pi\left( \bY^{(h)}\right)=(0,\ldots,h)\right)=\Prob\left(\Pi\left( \bX^{(h)}\right)=\Pi\left( \bY^{(h)}\right)=(h,\ldots,0)\right).
\end{align*}
%Finally, we obtain
%\begin{align*}
%\OPD_h\left( \bX^{(h)} , \bY^{(h)} \right)=&\frac{2\tau\left(\tilde{X}_{1},\ldots,\tilde{X}_{h},\tilde{Y}_{1},\ldots,\tilde{Y}_{h}\right)\sqrt{ \tilde{p}_{\tilde{\bX}^{(h)}} \left(1-\tilde{p}_{\tilde{\bX}^{(h)}}\right)\tilde{p}_{\tilde{\bY}^{(h)}}\left(1-\tilde{p}_{\tilde{\bY}^{(h)}}\right)}}{1-\sum_{\pi\in S_h}\Prob\left(\Pi\left( \bX^{(h)}\right)=\pi\right)\Prob\left(\Pi\left( \bY^{(h)}\right)=\pi\right)}\\
%&+\frac{\sum_{\pi\in S_h\setminus T_h} \Prob\left(\Pi\left( \bX^{(h)}\right)=\Pi\left( \bY^{(h)}\right)=\pi\right)}{1-\sum_{\pi\in S_h}\Prob\left(\Pi\left( \bX^{(h)}\right)=\pi\right)\Prob\left(\Pi\left( \bY^{(h)}\right)=\pi\right)}\\
%&-\frac{\sum_{\pi\in S_h\setminus T_h} \Prob\left(\Pi\left( \bX^{(h)}\right)=\pi\right)\Prob\left(\Pi\left( \bY^{(h)}\right)=\pi\right)}{1-\sum_{\pi\in S_h}\Prob\left(\Pi\left( \bX^{(h)}\right)=\pi\right)\Prob\left(\Pi\left( \bY^{(h)}\right)=\pi\right)}.
%\end{align*}
Now let $\pi=\left(\pi_0,\ldots,\pi_{h}\right)$ be a permutation of $0, \cdots, h$. If $\Pi\left(X_1,\ldots,X_{1+h}\right)=\pi$, it holds that $\left\{X_{1+\pi_0}\geq X_{1+\pi_1}\geq\ldots\geq X_{1+\pi_{h}}\right\}$. As a result, ordinal pattern dependence can be expressed by the following formula:
\begin{align*}
\OPD_h\left( \bX^{(h)} , \bY^{(h)} \right)
=&\frac{\sum_{\pi\in S_h} \Prob\left(\Pi\left(X_1,\ldots,X_{1+h}\right)=\Pi\left(Y_1,\ldots,Y_{1+h}\right)=\pi\right)- p_{\bX^{(h)},\pi}p_{\bY^{(h)},\pi}}{1-  \sum_{\pi\in S_h} \Prob\left(\Pi\left(X_1,\ldots,X_{1+h}\right)=\pi\right) \Prob\left(\Pi\left(Y_1,\ldots,Y_{1+h}\right)=\pi\right)}\\
=& \frac{ \sum_{\pi\in S_h} \Prob\left( X_{1+\pi_0}\geq X_{1+\pi_1}\geq\ldots\geq X_{1+\pi_{h}},Y_{1+\pi_0}\geq Y_{1+\pi_1}\geq\ldots\geq Y_{1+\pi_{h}}  \right)-p_{\bX^{(h)},\pi}p_{\bY^{(h)},\pi} }{1-\sum_{\pi\in S_h} p_{\bX^{(h)},\pi}p_{\bY^{(h)},\pi}}\\
=&\frac{ \sum_{\pi\in S_h} \Prob\left( X_{1+\pi_1}- X_{1+\pi_0}\leq 0, \ldots,X_{1+\pi_{h}}-X_{1+\pi_{h-1}}\leq 0, \ldots,Y_{1+\pi_{h}}-Y_{1+\pi_{h-1}}\leq 0  \right)-p_{\bX^{(h)},\pi}p_{\bY^{(h)},\pi} }{1-\sum_{\pi\in S_h} p_{\bX^{(h)},\pi}p_{\bY^{(h)},\pi}}\\
=&\frac{ \sum_{\pi \in S_h} \tau\left(  X_{1+\pi_1}- X_{1+\pi_0},\ldots, Y_{1+\pi_{h}}-Y_{1+\pi_{h-1}} \right)\sqrt{p_{\bX^{(h)},\pi}\left(1-p_{\bX^{(h)},\pi}\right)p_{\bY^{(h)},\pi}\left(1-p_{\bY^{(h)},\pi}\right)}}{1-\sum_{\pi\in S_h} p_{\bX^{(h)},\pi}p_{\bY^{(h)},\pi}}.
\end{align*}
%with 
%\begin{align*}
%p_{\bX^{(h)},\pi}&=\Prob\left(\Pi\left(X_1,\ldots,X_{1+h}\right)=\pi\right)\\
%&=\Prob\left( X_{1+\pi_1}\geq X_{1+\pi_2}\geq\ldots\geq X_{1+\pi_{1+h}}\right)\\
%&=  \Prob\left( X_{1+\pi_2}- X_{1+\pi_1}\leq 0, \ldots,X_{1+\pi_{1+h}}-X_{1+\pi_{h}}\leq 0\right)
%\end{align*}
% and 
% \begin{align*}
%p_{\bY^{(h)},\pi}&=\Prob\left(\Pi\left(Y_1,\ldots,Y_{1+h}\right)=\pi\right)\\
%&=\Prob\left( Y_{1+\pi_1}\geq Y_{1+\pi_2}\geq\ldots\geq Y_{1+\pi_{1+h}}\right)\\
%&=  \Prob\left( Y_{1+\pi_2}- Y_{1+\pi_1}\leq 0, \ldots,Y_{1+\pi_{1+h}}-Y_{1+\pi_{h}}\leq 0\right).
%\end{align*}
\end{proof}

\section{Conclusion and Outlook}
We have shown that ordinal pattern dependence is a multivariate measure of dependence in an axiomatic sense. 
When applied to bivariate time series,
it
can be interpreted as  a value describing
the co-movement of the  two component time series in equal moving windows.
In contrast to other dependence measures, it has thus been
developed against a time series background.  To make ordinal pattern dependence comparable to other multivariate dependence measures, we adapted the definitions of the latter to the time series approach. We figured out that univariate dependence measures do not carry enough information for an analysis of dependencies between two random vectors. They do not take any dynamical dependence into account, given for example by the cross-correlations of the considered random vectors. The same holds true for the multivariate extension of Pearson's $\rho$. Hence, both measures are inappropriate in a time series context. For Gaussian observations, there is a close relationship between ordinal pattern dependence and the multivariate version of Kendall's $\tau$. We proved that for $h=1$ ordinal pattern dependence of two random vectors can be represented as Kendall's $\tau$ of the corresponding increment vectors. The provided simulations show that the values of Kendall's $\tau$ and ordinal pattern dependence of the same data set differ in concrete situations. This emphasizes that the two measures operate on different levels, i.e., one operates on the level of the original process, the other on the level of increments. Moreover, for this reason ordinal pattern dependence is insensitive with respect to time shifts: it only relies on the dependence between the corresponding increments.
 However, an extension of ordinal pattern dependence that is sensitive to dependence shifted in time is given in \cite{schnurr:dehling:2017}. The authors of that article introduced time shifted ordinal pattern dependence that allows for time shifts between the moving windows of the two time series. Using this approach it is possible to detect a co-movement of the two time series that does not happen in the same moving window. 
Furthermore, it is an interesting topic for further research to compare ordinal pattern dependence to copula-based dependence measures. Originally, copulas were introduced to focus on the dependence \textit{within} a multivariate random vector without taking the marginal distributions into account. If all univariate margins are continuous, the multivariate generalizations of Kendall's $\tau$ and multivariate Spearman's $\rho$ in \cite{grothe:2014} only depend on the underlying copula. As these two measures are defined to measure dependence between two random vectors, this approach is a promising starting point to extend these ideas to a time series background. An investigation of ordinal pattern dependence with respect to the framework of copulas as well as a comparison to multivariate Spearman's $\rho$ is ongoing research, but not in the scope of the present article. 
Since ordinal pattern dependence admits a canonical interpretation, and since limit theorems in the short and the long-range dependent framework are at hand, we suggest to use this measure to complement the classical time series analysis with an ordinal point of view.
 
\section*{Acknowledgments} 
 We would like to thank two anonymous referees whose comments have helped to improve the presentation of our results. 
This research was supported in part by the German Research Foundation (DFG) through  Collaborative Research Center SFB 823   Statistical Modelling of Nonlinear Dynamic Processes and the project Ordinal-Pattern-Dependence: Grenzwerts\"atze und Strukturbr\"uche im langzeitabh\"angigen Fall mit Anwendungen in Hydrologie, Medizin und Finanzmathematik (SCHN 1231/3-2).

\bibliographystyle{myjmva} 
%%  \bibliography{<your bibdatabase>}
%\newpage
\bibliography{DepMeasuresRevision_July2021arXiv}
%\end{thebibliography}

\appendix
\newpage

\newgeometry{paper=a4paper, left=30mm, right=0mm, top=50mm, bottom=50mm}

\begin{landscape}
\section{Simulation study}
In this section, we empirically compare multivariate Kendall's $\tau$  and ordinal pattern dependence by simulation studies for different settings. 

\begin{table}[ht]
\centering
\footnotesize
\caption{We  generate 	$X_i$, $1\leq i\leq n$, and 	$Y_i$, $1\leq i\leq n$, as two independent AR(1) time series with
$X_i=\rho X_{i-1}+\varepsilon_i, \ Y_i=\rho Y_{i-1}+\eta_i,$
where $\left|\rho\right|<1$, and $\varepsilon_i$,  $1\leq i\leq n$,   and $\eta_i$,  $1\leq i\leq n$,  are two independent  sequences of random variables, both i.i.d. standard normally distributed.
}
\hrule
\vspace{1em}
\begin{tabular}{llrrrrrrrrrrrrrrrrrrrrrr}
& & & $\rho=0.1$ &  & & $\rho=0.5$ & & & $\rho=0.9$ \\ 
				\cline{4-5}   \cline{7-8} \cline{10-11} \\
  method & $n$ & $h$ & mean (sd) & median (IQR) & & mean (sd) & median (IQR) & & mean (sd) & median (IQR) \\ 
  \hline
  \\
  \multirow{3}{*}{\rotatebox{90}{\footnotesize Kendall's $\tau$}} 
 & $100$ & 1 & 0.007  (0.062) & 0.004 (0.086) & & 0.012 (0.087) & 0.011 (0.121) && 0.018 (0.168) & 0.017 (0.238) \\ 
   & $300$ & 1 & 0.003 (0.035) & 0.003 (0.048) & & 0.004 (0.051) & 0.004 (0.068) & & 0.008 (0.111) & 0.007 (0.153) \\ 
   & $500$ & 1 & 0.002 (0.027) & 0.002 (0.036) && 0.003 (0.040) & 0.002 (0.053) && 0.006 (0.088) & 0.004 (0.119) \\ 
   \\
   \multirow{3}{*}{\rotatebox{90}{\footnotesize opd}}   & $100$ & 1 & 0.009 (0.106) & 0.003 (0.139)& & 0.011 (0.101) & 0.016 (0.137) && 0.008 (0.098) & 0.003 (0.137) \\ 
   & $300$  & 1 & 0.003 (0.062) & 0.002  (0.085) & & 0.004 (0.059) & 0.005 (0.083) & & 0.003 (0.058) & 0.004 (0.075)\\ 
   & $500$  & 1 & 0.002 (0.048) & 0.004 (0.062) && 0.003 (0.046) & 0.003 (0.063)& & 0.001 (0.044) & -0.000 (0.061) \\ 
       \\
    \hline
    \\
   \multirow{3}{*}{\rotatebox{90}{\footnotesize Kendall's $\tau$}}   & $100$ & 2 & 0.008 (0.053) & 0.002 (0.070) & & 0.015 (0.083) & 0.011 (0.114) & & 0.031 (0.170) & 0.027 (0.239) \\ 
   & $300$  & 2 & 0.003 (0.030) & 0.002 (0.040) & & 0.005 (0.048) & 0.004 (0.065) & & 0.013 (0.111) & 0.011 (0.154) \\ 
   & $500$  & 2 & 0.002 (0.023) & 0.001 (0.031) & & 0.003 (0.038) & 0.002 (0.050) & & 0.009 (0.089) & 0.007  (0.120) \\
   \\ 
   \multirow{3}{*}{\rotatebox{90}{\footnotesize opd}}   & $100$ & 2 & 0.003 (0.055) & -0.000 (0.073) & & 0.004 (0.055) & 0.001 (0.077) && 0.004 (0.056) & 0.002 (0.075) \\ 
   & $300$  & 2 & 0.001 (0.033) & -0.000 (0.044)& & 0.002 (0.032) & 0.001 (0.044) & & 0.001 (0.033) & 0.001 (0.044) \\ 
   & $500$  & 2 & 0.000 (0.025) & 0.000 (0.034) & & 0.001 (0.025) & 0.001 (0.034) & & 0.000 (0.025) & -0.001 (0.034) \\ 
       \\
    \hline
    \\
  \multirow{3}{*}{\rotatebox{90}{\footnotesize Kendall's $\tau$}}    & $100$ & 3 & 0.006 (0.043) & -0.002 (0.053) && 0.015 (0.076) & 0.007 (0.103) && 0.040 (0.169) & 0.031 (0.239) \\ 
   & $300$  & 3 & 0.002 (0.024) & -0.000 (0.032) & & 0.005 (0.044) & 0.002 (0.059) && 0.016 (0.111) & 0.014 (0.154) \\ 
   & $500$  & 3 & 0.001 (0.019) & 0.000 (0.024) & & 0.003 (0.035) & 0.001 (0.047) && 0.012 (0.089) & 0.009 (0.120) \\ 
   \\
    \multirow{3}{*}{\rotatebox{90}{\footnotesize opd}}  & $100$ & 3 & 0.001 (0.027) & -0.001 (0.035) & & 0.001 (0.027) & -0.002 (0.037) & & 0.002 (0.031) & -0.001 (0.041) \\ 
   & $300$  & 3 & 0.001 (0.015) & -0.001 (0.021) & & 0.001 (0.016) & -0.000 (0.021) & & 0.001 (0.018) & -0.001 (0.024) \\ 
    & $500$  & 3 & -0.000 (0.012) & -0.001 (0.016) & & 0.000 (0.013) & -0.000 (0.017)& & 0.000  (0.014) & -0.000 (0.019) \\ \\
    \hline
\end{tabular}
\label{table:independent}
\end{table}
\end{landscape}

\begin{landscape}
\begin{table}[ht]
\centering
\footnotesize
\caption{ We simulate $(X_i)$, $i\geq 1$, and $(Y_i)$, $i\geq 1$, as  sequences of independent,  multivariate normal random vectors with values in $\mathbb{R}^3$ and a joint normal distribution with expectation $0$ and covariance matrix \eqref{eq:cov}.
%$\Sigma$ is positive definite if $\rho\in\{0.1, 0.2, 0,3\}$.
}\label{table:mult_normal}
\hrule
\vspace{1em}
\begin{tabular}{llrrrrrrrrrrrrrrrrrrrrrr}
& & & $\rho=0.1$ &  & & $\rho=0.2$ & & & $\rho=0.3$ \\ 
				\cline{4-5}   \cline{7-8} \cline{10-11} \\
  method & $n$ & $h$ & mean (sd) & median (IQR) & & mean (sd) & median (IQR) & & mean (sd) & median (IQR) \\ 
  \hline
  \\
  \multirow{3}{*}{\rotatebox{90}{\footnotesize Kendall's $\tau$}}  & $100$  & 1 & 0.045 (0.034) & 0.045 (0.047) & & 0.091 (0.037) & 0.091 (0.050) & & 0.141 (0.039) & 0.141 (0.053) \\ 
   & $300$  & 1 & 0.045 (0.020) & 0.045 (0.026) & & 0.091 (0.021) & 0.091 (0.028) & & 0.142 (0.022) & 0.142 (0.030) \\ 
   & $500$  &  1 & 0.045 (0.015) & 0.044 (0.021) & & 0.091 (0.016) & 0.091 (0.023) & & 0.142 (0.017) & 0.142 (0.023) \\
   \\ 
    \multirow{3}{*}{\rotatebox{90}{\footnotesize opd}}     & $100$   &  1 & 0.066 (0.064) & 0.066 (0.081) && 0.129 (0.063) & 0.130 (0.086) & & 0.196 (0.066) & 0.198 (0.087) \\ 
   & $300$ &  1 & 0.065 (0.036) & 0.064 (0.048) & & 0.128 (0.037) & 0.129 (0.051) & & 0.195 (0.038) & 0.195 (0.053) \\ 
   &  $500$   &  1 & 0.065 (0.029) & 0.065 (0.038) & & 0.128 (0.028) & 0.128 (0.039) & & 0.195 (0.030) & 0.194 (0.040) \\ 
   \\
      \hline
   \\
  \multirow{3}{*}{\rotatebox{90}{\footnotesize Kendall's $\tau$}}    & $100$   &  2 & 0.031  (0.030) & 0.029 (0.041) && 0.063 (0.033) & 0.062  (0.045) && 0.100 (0.036) & 0.100 (0.050) \\ 
   & $300$  &  2 & 0.030 (0.017) & 0.029 (0.023) & & 0.062 (0.019) & 0.062 (0.026) & & 0.100 (0.021) & 0.100 (0.029) \\ 
   &  $500$   &  2 & 0.030 (0.013) & 0.030 (0.018) & & 0.063 (0.015) & 0.062 (0.019) & & 0.100 (0.016) & 0.100 (0.022) \\ 
   \\
    \multirow{3}{*}{\rotatebox{90}{\footnotesize opd}}     & $100$   &  2 & 0.031 (0.034) & 0.031  (0.047) & & 0.065 (0.037) & 0.063 (0.051) & & 0.101  (0.040) & 0.099 (0.054) \\ 
   & $300$  &  2 & 0.030 (0.020) & 0.030 (0.027) && 0.064 (0.022) &  0.063 (0.029) & & 0.102 (0.024) & 0.101 (0.031) \\ 
   &  $500$   &  2 & 0.030 (0.016) & 0.030 (0.021) & & 0.064 (0.017) & 0.063 (0.023) & &  0.102 (0.018) & 0.102 (0.025) \\ 
   \\
      \hline
   \\
   \multirow{3}{*}{\rotatebox{90}{\footnotesize Kendall's $\tau$}}   & $100$   &  3 & 0.019 (0.024) & 0.017 (0.032) & & 0.041 (0.029) & 0.039 (0.038) && 0.067 (0.033) & 0.064 (0.043) \\ 
   & $300$ &  3 & 0.019 (0.014)& 0.018 (0.019) && 0.042 (0.017) & 0.041 (0.022) && 0.068 (0.019) & 0.067 (0.026) \\ 
   &  $500$   &  3 & 0.019 (0.011) & 0.018 (0.015) && 0.041 (0.013) & 0.041 (0.018) && 0.069 (0.015) & 0.068 (0.020) \\ 
   \\
    \multirow{3}{*}{\rotatebox{90}{\footnotesize opd}}     & $100$   &  3 & 0.011 (0.017) & 0.010  (0.023) && 0.024 (0.020) & 0.023 (0.028) && 0.041  (0.023) & 0.040 (0.032) \\ 
   & $300$  &  3 & 0.011 (0.010) & 0.011 (0.014) && 0.024 (0.012) & 0.024 (0.015) && 0.041 (0.013) & 0.040 (0.018) \\ 
   &  $500$   &  3 & 0.011 (0.008) & 0.011 (0.010) && 0.024 (0.009) & 0.024 (0.012) && 0.041 (0.010) & 0.041 (0.014) \\ \\
   \hline
\end{tabular}
\end{table}
\end{landscape}

\begin{landscape}
\begin{table}[ht]
\centering
\footnotesize
	\caption{We generate
			$X_i$, $1\leq i\leq n$, from an AR(1) time series with parameter $\rho$, while $Y_i$, $1\leq i\leq n$, corresponds to $X_i$, $1\leq i\leq n$, shifted by one time point, i.e.,
 $Y_i=X_{i+1}$
.}\label{table:shift}
\hrule 
\vspace{1em}
\begin{tabular}{llrrrrrrrrrrrrrrrrrrrrrr}
				& & & $\rho=0.1$ &  & & $\rho=0.5$ & & & $\rho=0.9$ \\ 
				\cline{4-5}   \cline{7-8} \cline{10-11} \\
  method & $n$ & $h$ & mean (sd) & median (IQR) & & mean (sd) & median (IQR) & & mean (sd) & median (IQR) \\ 
  \hline
  \\
  \multirow{3}{*}{\rotatebox{90}{\footnotesize Kendall's $\tau$}} &
$100$ &  1 & 0.354 (0.059) & 0.354 (0.080) & & 0.509 (0.058) & 0.512 (0.077) & & 0.738  (0.055) & 0.744 (0.072) \\ 
  & $300$ &   1 & 0.360 (0.034) & 0.361 (0.046) & & 0.521 (0.033) & 0.522 (0.044) & & 0.777  (0.030) & 0.779 (0.039) \\ 
  & $500$ &  1 & 0.361 (0.026) & 0.361 (0.036) & & 0.524 (0.025) & 0.525 (0.033) & & 0.784  (0.023) & 0.786 (0.031) \\ 
  \\
   \multirow{3}{*}{\rotatebox{90}{\footnotesize opd}} &  $100$ &  1 & -0.282 (0.083) & -0.279 (0.112) & & -0.152 (0.090) & -0.157 (0.118) & & -0.026 (0.097) & -0.023  (0.131) \\ 
 &   $300$  &   1 & -0.292 (0.049) & -0.292 (0.066) & & -0.159 (0.053) & -0.159 (0.071) & & -0.030 (0.057) & -0.028 (0.074) \\ 
  &  $500$  &   1 & -0.295 (0.038) & -0.295 (0.050) & & -0.158 (0.041) & -0.157 (0.054) & & -0.030 (0.044) & -0.029 (0.060) \\
    \\
    \hline
    \\
    \multirow{3}{*}{\rotatebox{90}{\footnotesize Kendall's $\tau$}}   &  $100$ &   2 & 0.437 (0.070) & 0.438 (0.095) & & 0.575 (0.065) & 0.579 (0.086) && 0.783 (0.055) & 0.791 (0.071) 
      \\
  &  $300$  &  2 & 0.450 (0.038) & 0.450 (0.052) & & 0.591 (0.036) & 0.593 (0.049) && 0.816 (0.027) & 0.817 (0.036) \\ 
  &  $500$  &   2 & 0.453 (0.030) & 0.454 (0.041) & & 0.594 (0.027) & 0.595 (0.037) && 0.823 (0.021) & 0.824 (0.028) \\
  \\ 
     \multirow{3}{*}{\rotatebox{90}{\footnotesize opd}} &  $100$ &   2 & -0.088 (0.037) & -0.088 (0.049) & & -0.030 (0.044) & -0.031 (0.060) && 0.050  (0.052) & 0.049 (0.070) \\ 
   & $300$  &  2 & -0.088 (0.021) & -0.088 (0.028) & & -0.028 (0.025) & -0.029 (0.033) && 0.052 (0.030) & 0.051 (0.041) \\ 
   & $500$  &  2 & -0.087 (0.016) & -0.087 (0.022) & & -0.028 (0.019) & -0.028 (0.025) && 0.052 (0.024) & 0.052 (0.032) \\ 
    \\
    \hline
    \\
      \multirow{3}{*}{\rotatebox{90}{
 \footnotesize Kendall's $\tau$}} &   $100$ &  3 & 0.462 (0.083) & 0.465 (0.113) & & 0.602 (0.075) & 0.609 (0.101) && 0.804 (0.058) & 0.814 (0.073) \\ 
  &  $300$  &   3 & 0.485 (0.047) & 0.485 (0.062) & & 0.624 (0.040) & 0.625 (0.053) && 0.837 (0.027) & 0.840 (0.035) \\ 
  &  $500$  &   3 & 0.489 (0.035) & 0.489 (0.047) & & 0.629 (0.030) & 0.630 (0.040) && 0.844 (0.020) & 0.845 (0.026) \\ 
  \\
   \multirow{3}{*}{\rotatebox{90}{
 \footnotesize  opd}} &   $100$ &  3 & -0.029 (0.016) & -0.030 (0.022) & & -0.008 (0.024) & -0.010 (0.033) && 0.040 (0.038) & 0.037 (0.050) \\ 
  &  $300$  &  3 & -0.025 (0.010) & -0.025 (0.013) & & -0.003 (0.014) & -0.004 (0.019) && 0.045 (0.022) & 0.045 (0.030) \\ 
   & $500$  &   3 & -0.024 (0.007) & -0.024 (0.010) & & -0.002 (0.011) & -0.002 (0.015) && 0.046 (0.017) & 0.046 (0.023) \\ \\
   \hline
\end{tabular}
\end{table}
\end{landscape}
\end{document}